\documentclass[onefignum,onetabnum,preprint]{siamart190516}

\usepackage{lipsum}
\usepackage{amsfonts}
\usepackage{graphicx}
\usepackage{epstopdf}
\usepackage{algorithmic}
\ifpdf
  \DeclareGraphicsExtensions{.eps,.pdf,.png,.jpg}
\else
  \DeclareGraphicsExtensions{.eps}
\fi

\newsiamremark{remark}{Remark}
\newsiamremark{hypothesis}{Hypothesis}
\crefname{hypothesis}{Hypothesis}{Hypotheses}
\newsiamthm{claim}{Claim}

\headers{Multi-fidelity surrogate modeling for time-series outputs}{B. Kerleguer}

\title{Multi-fidelity surrogate modeling for time-series outputs }
\author{Baptiste Kerleguer \footnotemark[2]\ \footnotemark[3] }

\usepackage{amsopn}
\DeclareMathOperator{\diag}{diag}

\ifpdf
\hypersetup{
  pdftitle={Multi-fidelity surrogate modeling for time-series outputs},
  pdfauthor={B. Kerleguer}
}
\fi

\usepackage[utf8]{inputenc}
\usepackage[T1]{fontenc}		
\usepackage[english]{babel}

\usepackage{enumitem}
\usepackage{bm}

\newcommand{\R}{\mathbb{R}}

\newcommand{\E}{\mathbb{E}}
\newcommand{\V}{\mathbb{V}}
\newcommand{\Cov}{\mathrm{Cov}}

\begin{document}

\maketitle

\renewcommand{\thefootnote}{\fnsymbol{footnote}}

\footnotetext[2]{CEA,DAM,DIF, F-91297, Arpajon, France \email{baptiste.kerleguer@cea.fr}}
\footnotetext[3]{Centre de Math\'ematiques Appliqu\'ees, Ecole Polytechnique, Institut Polytechnique de Paris, 91128 Palaiseau Cedex, France}

\renewcommand{\thefootnote}{\arabic{footnote}}

\begin{abstract}
  This paper considers the surrogate modeling of a complex numerical code in a multi-fidelity framework when the code output is a time series \textcolor{black}{and two code levels are available: a high-fidelity and expensive code level and a low-fidelity and cheap code level.}
\textcolor{black}{The goal is to emulate a fast-running approximation of the high-fidelity code level.}
\textcolor{black}{A}n original Gaussian process regression method is proposed \textcolor{black}{that uses} an experimental design of the low- and high-fidelity code levels.
The code output is expanded on a basis built from the experimental design.
The first coefficients of the expansion of the code output are processed by a co-kriging approach.
The last coefficients are processed by a kriging approach with covariance tensorization.
The resulting surrogate model \textcolor{black}{provides a predictive mean and a predictive variance of the output of the high-fidelity code level.}
\textcolor{black}{It} is shown to have better performance in terms of prediction errors and uncertainty quantification than standard dimension reduction techniques.
\end{abstract}



\section{Introduction}

Advances in scientific modeling have led to the development of complex and computationally expensive codes.
\textcolor{black}{To solve the problem of computation time for tasks such as optimization or calibration of complex numerical codes, surrogate models are used.}
The surrogate modeling approach consists in building a surrogate model of a complex numerical code from a data set computed from an experimental design.
A well-known method to build surrogate models is Gaussian process \textcolor{black}{(GP)} regression.
This method, also called kriging, was originally proposed by \cite{krige1951statistical} for geostatistics.
This method has subsequently been used for computer experiments and in particular in the field of uncertainty quantification \textcolor{black}{(UQ)}, see \cite{santner2003design,williams2006gaussian}.

\textcolor{black}{It is common for complex codes to have different versions that are more or less accurate and more or less computationally expansive.}
The particular case that interests us is when codes are hierarchical, i.e. they are classified according to their computational cost and their accuracy.
The more accurate the code, the more expensive it is.
The autoregressive scheme presented by \cite{kennedy2000predicting} is the first \textcolor{black}{major} result in the field of \textcolor{black}{multi-fidelity} Gaussian process regression.
This technique has been amended by \cite{le2014recursive} in order to reduce the overall co-kriging problem  to several independent kriging problems.
The papers \cite{goh2013prediction, pilania2017multi} present different application cases and \cite{giselle2019issues} is a synthesis of the use of multi-fidelity for surrogate modeling.
\textcolor{black}{In \cite{ma2020objective} the author introduces objective prior for the hyperparameters of the autoregressive model.}

\textcolor{black}{New methods for multi-fidelity surrogate modeling have been introduced.
\textcolor{black}{In particular,} Deep Gaussian processes have been proposed to solve cases where the interactions between code levels are more complex \cite{perdikaris2017nonlinear}.
This type of methods \textcolor{black}{can deal with UQ} \cite{cutajar2019deep} but thus are time-consuming and do not scale up the dimension of outputs.
Neural networks have also been used to emulate multi-fidelity computer codes.
In particular, the method proposed in \cite{meng2020composite} is a neural network method with an AR(1)-based model and is scalable to high-dimensional outputs.
However, UQ is not taken into account in this method.}

Among the codes with high-dimensional outputs, we are interested in those whose outputs are functions of one variable.
When they are sampled, such outputs are called time series.
Previous work has solved the problem of functional outputs only in the \textcolor{black}{single}-fidelity case.
Two methods have been considered to solve the \textcolor{black}{single}-fidelity problem: reduce the dimension of the outputs \cite{nanty2017uncertainty} or adapt the \textcolor{black}{covariance} kernel \cite{rougier2008efficient}.
\textcolor{black}{In the context} of dimension reduction\textcolor{black}{, surrogate models generally neglect the uncertainty quantification in the basis, which can be problematic for the quantification of prediction uncertainty}.
Moreover, large data sets (containing many low-fidelity data) lead to ill-conditioned covariance matrices that are difficult to invert.
As proposed in \cite{perrin2019adaptive}, it is possible to strongly constrain the covariance kernel which makes it possible to improve the estimation compared to the dimension reduction method.
However, this method implies that the covariance must be separable, which reduces the use cases.
Knowing that the AR(1) multi-fidelity model for \textcolor{black}{GP} regression uses co-kriging, \cite{nerini2010cokriging} presents an interesting approach for co-kriging in the context of functional outputs, which is based on dimension reduction.
An approach to multi-fidelity with functional outputs is presented in \cite{grujic2018cokriging} for multivariate Hilbert space valued random fields.

In this work, we introduce an original approach to the construction of a surrogate model in the framework of hierarchical multi-fidelity codes with time-series outputs.
The main idea is to combine a reduction method of the output dimension, that fits well with the autoregressive model of multi-fidelity co-kriging, and a \textcolor{black}{special single-}fidelity method that allows to treat time series output by \textcolor{black}{GP} regression  with covariance tensorization.
\textcolor{black}{We address the case of one high-fidelity code and one low-fidelity code.}

In \Cref{sec::Regressiontools} we give the main elements of the \textcolor{black}{GP} regression theory that are needed in our paper.

In \Cref{sec::AR1multiproj} we develop a \textcolor{black}{Cross-Validation Based (CVB)} method \textcolor{black}{ to determine an appropriate basis} based on K-fold cross-validation which uses only low-fidelity data and allows us to \textcolor{black}{study} the first two moments of the basis vectors. 
Once projected onto the reduced space it is possible to characterize the models of the first coefficients of the expansion of the code output onto the basis by the multi-fidelity \textcolor{black}{GP} regression method.

\textcolor{black}{The projection of the code output onto the orthogonal of the reduced space forms the orthogonal part.}
\textcolor{black}{If we neglect the orthogonal part, then we obtain the results of \Cref{sec::AR1multiproj}.
However, the orthogonal part may not be negligible.
In \Cref{sec::AR1multiTens} we treat by GP} regression with covariance  tensorization \textcolor{black}{the orthogonal part}.
The latter approach collectively addresses the orthogonal part of the high-fidelity code\textcolor{black}{,} and it makes it possible  to better predict the output of the high-fidelity code and to better quantify the \textcolor{black}{prediction uncertainty}.

The results presented in the numerical example in \Cref{sec::Illustration} confirm that the processing of the orthogonal part is important.
In this example we test the different methods presented in the paper \textcolor{black}{as well as NN state of the art multi-fideliy methods,} and we assess their performance in terms of prediction errors and uncertainty quantification.


\section{Gaussian Process regression}
\label{sec::Regressiontools}
\textcolor{black}{The autoregressive multi-fidelity model has been introduced by \cite{kennedy2000predicting}.
The authors in \cite{le2014recursive} have simplified the computation.
We present our AR model in \cref{subsec::multiGPreg}.
In parallel, GP regression has been used with covariance tensorization in \cite{rougier2008efficient} and improved in \cite{conti2009gaussian, perrin2019adaptive}.
This method allows to extend GP regression to time-series outputs.
We have exploited this method to build our own methodology for time-series regression and we present it in \cref{subsec::tensRegcov}}

\textcolor{black}{Let us consider a complex numerical code where the input is a point $\bm{x} \in Q$, with $Q$ being a domain in $\R^d$ and the output is a function of a one-dimensional variable. 
We are interested in hierarchical codes, which means that there are several code levels that can be classified according to their fidelity.
In this work, we focus on only two code levels, a high-fidelity code and a low-fidelity code.
In what follows, $\text{H}$ represents the high-fidelity code and $\text{L}$ the low-fidelity code, the generic notation is $\text{F} \in \{ \text{H},\text{L}\}$.
For any given  input $\bm{x}$, we can run the $\text{F}$ code and observe the data $z_\text{F}( \bm{x}, t)$  for all $t$ in a fixed regular grid $\left\{t_u, u=1, \ldots, N_t\right\}$ in $[0,1]$.
However, the cost of the code allows only a limited number of code calls.
This induces the use of the experimental design $D_\text{F} = \left\{\bm{x}^{(1)}, \ldots, \bm{x}^{(N_\text{F})}\right\}$.
$N_\text{F}$ is the number of observations of the code $\text{F}$.
The $N_t\times N_\text{F}$ matrix containing the observations for $\bm{x} \in D_\text{F}$ is $\bm{Z}_{\text{obs}}^\text{F}$.
Our goal is to predict the values of $\left( z_\text{H}(\bm{x},t_u)\right)_{u=1,\ldots,N_t}$ given $\left(\bm{Z}_{\text{obs}}^\text{H}, \bm{Z}_{\text{obs}}^\text{L}\right)$ for a point $\bm{x}\in Q$ with the quantification of the prediction uncertainty.
We model the prior knowledge of the code output $\left(z_\text{L},z_\text{H}\right)$ as a Gaussian process $\left(Z_\text{L},Z_\text{H}\right)$.
We denote by $\mathcal{Z}^{\text{F}}$ the random vector containing the random variables $\left(Z_\text{F}(\bm{x}^{(i)},t_u)\right)_{\substack{u=1,\ldots,N_t\\ u=1,\ldots,N_t}}$.
The combination of $\mathcal{Z}^{\text{L}}$ and $\mathcal{Z}^{\text{H}}$ is $\mathcal{Z}$.
}

\subsection{Multi-fidelity Gaussian process regression} 
\label{subsec::multiGPreg}
In this section we want to build a surrogate model of a code $\alpha_\text{H}(\bm{x})$ whose input $\bm{x}$ is in $Q \subset \R^d$ and whose scalar output is in $\R$.
The construction of a surrogate model for complex computer code is difficult because of the lack of available experimental outputs.
We consider the situation in which a cheaper and approximate code $\alpha_\text{L}(\bm{x})$ is available.
In this section, we apply the regression method presented by \cite{kennedy2000predicting}, reviewed in \cite{forrester2007multi} and improved in \cite{le2014recursive}.

\textcolor{black}{We model the prior knowledge of t}he code output $(\alpha_\text{L},\alpha_\text{H})$ \textcolor{black}{as a} Gaussian process $\left(A_\text{L}, A_\text{H}\right)$.
The vector containing the values of $\alpha_\text{F}(\bm{x})$ at the points of the experimental design $D_\text{F}$ are denoted by $\alpha^{\text{F}}$
and $\mathcal{A}^{\text{F}}$ is the Gaussian vector containing \textcolor{black}{ $A_\text{F}(\bm{x})$, $\bm{x}\in D_\text{F}$.}
The combination of $\mathcal{A}^{\text{L}}$ and $\mathcal{A}^{\text{H}}$ is $\mathcal{A}$.
So is $\alpha$, the combination of $\alpha^{\text{L}}$ and $\alpha^{\text{H}}$.
We present the recursive model of multi-fidelity introduced by \cite{le2014recursive}.
The experimental design is constructed such that $D_\text{H} \subset D_\text{L}$.
We assume the low-fidelity code is computationally cheap, and that we have access to a large experimental design \textcolor{black}{for the low-fidelity code}, i.e. $N_\text{L}\gg N_\text{H}$.


We consider the hierarchical model introduced by \cite{le2014recursive}:

\begin{equation}\label{multikrigeq}
\left\{
\begin{array}{rcl}
A_\text{H}(\bm{x}) & = & \rho_\text{L}(\bm{x})\tilde{A}_\text{L}(\bm{x})+\delta(\bm{x})\\ 
\tilde{A}_\text{L}(\bm{x}) & \perp & \delta(\bm{x})\\
\rho_\text{L}(\bm{x}) & = & g_\text{L}^T(\bm{x})\beta_\rho\\
\end{array}
\right. ,
\end{equation}
where $\perp$ means independence,
$~^T$ stands for the transpose,
\begin{equation}
    \label{eq::deltaGPmulti}
    \left[\delta(\bm{x}) | \beta_H,\sigma_H \right]\sim \mathcal{GP}\left(f^T_\text{H}(\bm{x})\beta_\text{H},\sigma_\text{H}^2r_\text{H}(\bm{x},\bm{x}^\prime)\right),
\end{equation}
and $\tilde{A}_\text{L}(\bm{x})$ is a Gaussian process conditioned by the values $\alpha^{\text{L}}$.
Its distribution is the one of $\left[A_\text{L}(\bm{x})|\mathcal{A}^{\text{L}}=\alpha^{\text{L}}, \beta_\text{L}, \sigma_\text{L}\right]$ with
\begin{equation}
    \left[  A_\text{L} (\bm{x}) | \beta_L,\sigma_L\right] \sim \mathcal{GP}\left(f^T_\text{L}(\bm{x})\beta_\text{L},\sigma_\text{L}^2r_\text{L}(\bm{x},\bm{x}^\prime)\right).
\end{equation}
Therefore, the distribution of $\tilde{A}_\text{L}(\bm{x})$ is Gaussian with mean $\mu_{\tilde{A}_\text{L}}(\bm{x})$ and variance $\sigma_{\tilde{A}_\text{L}}^2(\bm{x})$:
\begin{align}
    \mu_{\tilde{A}_\text{L}} (\bm{x})=&
    f_\text{L}^T(\bm{x})\beta_\text{L} + r_\text{L}^T (\bm{x}) C_\text{L}^{-1}\left(\alpha^{\text{L}}- F_\text{L}\beta_\text{L}\right),\\
    \sigma_{\tilde{A}_\text{L}}^2(\bm{x}) =&
    \sigma_\text{L}^2\big( 
    r_\text{L}(\bm{x},\bm{x}) - r_\text{L}^T (\bm{x}) C_\text{L}^{-1}r_\text{L} (\bm{x}) \big).
\end{align}
Here:
\begin{itemize}[noitemsep,topsep=0pt,parsep=0pt,partopsep=0pt]
    \item[-] $\mathcal{GP}$ means Gaussian process,
    \item[-] $g_\text{L}(\bm{x})$ is a vector of $q_\text{L}$ regression functions,
    \item[-] $f_\text{F}(\bm{x})$ are vectors of $p_\text{F}$ regression functions,
    \item[-] $r_\text{F}(\bm{x},\bm{x}^\prime)$ are correlation functions,
    \item[-] $\beta_\text{F}$ are $p_\text{F}$-dimensional vectors,
    \item[-] $\sigma_\text{F}^2$ are \textcolor{black}{positive real numbers},
    \item[-] $\beta_\rho$ is a $q$-dimensional vector of adjustment parameters,
    \item[-] $C_\text{F} = \left(r_\text{F}(\bm{x}^{(i)}, \bm{x}^{(j)})\right)_{i,j=1}^{N_\text{F}}$ is the $N_\text{F}\times N_\text{F}$ correlation matrix of $\mathcal{A}^{\text{F}}$,
    \item[-] $r_\text{F}(\bm{x}) = \left(r_\text{F}(\bm{x},\bm{x}^{(i)})\right)_{i=1}^{N_\text{F}}$ is the $N_\text{F}$-dimensional vector of correlations between $A_\text{F}(\bm{x})$ and $\mathcal{A}^{\text{F}}$,
    \item[-] $F_\text{F}$ is the $N_\text{F} \times p_\text{F}$ matrix containing the values of $f_\text{F}^T(\bm{x})$ for $\bm{x}\in D_{\text{F}}$.
\end{itemize}
For $\bm{x}\in Q$, the conditional distribution of $A_\text{H}(\bm{x})$ is:

\begin{equation}
\label{calautomodeleq}
\left[A_\text{H}(\bm{x}) |\mathcal{A}=\alpha,\beta,\beta_\rho,\sigma^2\right]\sim \mathcal{N}\left(\mu_{A_\text{H}}(\bm{x}),\sigma^2_{A_\text{H}}(\bm{x}) \right),
\end{equation}
where $\beta = \left(\beta_\text{H}^T, \beta_\text{L}^T\right)^T$ is the $p_\text{H}+p_\text{L}$-dimensional vector of regression parameters, $\sigma^2=\left(\sigma_\text{L}^2,\sigma_\text{H}^2\right)$ are the variance parameters,
\begin{align}
\nonumber
        \mu_{A_\text{H}}(\bm{x}) = & g_\text{L}^T(\bm{x})\beta_\rho \mu_{\tilde{A}_\text{L}}(\bm{x})
        + f_\text{H}^T(\bm{x})\beta_\text{H} \\
        &  + r_\text{H}^T (\bm{x}) C_\text{H}^{-1}\left(\alpha^{\text{H}}-\rho^\text{L}(D_\text{H}) \odot \alpha^{\text{L}}(D_\text{H})-F_\text{H}\beta_\text{H}\right)
\label{meanautomodeleq}
\end{align}
and
\begin{equation}\label{varautomodeleq}
\sigma^2_{A_\text{H}}(\bm{x}) =\left(g_\text{L}^T(\bm{x})\beta_\rho\right)^2 \sigma^2_{\tilde{A}_\text{L}}(\bm{x})
+ \sigma_\text{H}^2\left(1-r^T_\text{H}(\bm{x})C_\text{H}^{-1}r_\text{H}(\bm{x})\right).
\end{equation}
The notation $\odot$ is the element by element matrix product.
$\rho^\text{L}(D_\text{H})$ is the $N_\text{H}$-dimensional vector containing the values of $\rho_\text{L}(\bm{x})$ for $\bm{x}\in D_\text{H}$.
$\alpha^{\text{L}}\left(D_\text{H}\right)$ is the $N_\text{H}$-dimensional vector containing the values of $\alpha_\text{L}(\bm{x})$ at the points of $D_\text{H}$.

\textcolor{black}{The prior distributions of the parameters $\beta$ and $\sigma$ are given in \cref{asec::MultiFi}.}
The hyper-parameters of the covariance kernels $r_\text{L}$ and $r_\text{H}$ can be estimated by maximum likelihood or by leave-one-out cross validation \cite{bachoc2013cross}.
The nested property of the experimental design sets $D_\text{H}\subset D_\text{L}$ is not necessary to build the model but it is simpler to estimate the parameters with this assumption \textcolor{black}{\cite{zhou2020generalized}}.
Moreover, the ranking of codes and the low computer cost of the low-fidelity code allow for a nested design for practical applications.

\subsection{Gaussian process regression for functional outputs}
\label{subsec::tensRegcov}
In this subsection, we address \textcolor{black}{GP} regression for a simple-fidelity code with time-series output.
For the calculation of surrogate models with functional outputs, there are two different techniques.
The simplest ones are dimension reduction techniques as presented in \cite{aversano2019application,nanty2017uncertainty} (see \Cref{sec::AR1multiproj}).
An alternative is presented here, this method is \textcolor{black}{GP} regression with covariance tensorization.
The method is presented in \cite{rougier2008efficient} and the estimation of the hyper-parameters is from \cite{perrin2019adaptive}.

In this section \textcolor{black}{and the following ones} we consider that the output is a time-dependent function observed on a fixed time grid $\{ t_u\}_{u=1,\cdots,N_t}$, with $N_t\gg 1$, which is called a time series.

The experimental design in a times-series output case is very different from a scalar output case.
In particular, for a value $\bm{x}$ in the experimental design $D$, all the $t$ of the time grid are in the experimental design.
The $N_t \times N_x$ matrix containing the observations is $\bm{Z}_\text{obs} = \left(z(\bm{x}^{(i)},t_u)\right)_{\substack{u=1,\ldots,N_t\\ i=1,\ldots,N_x}}$.
In \textcolor{black}{GP} regression \textcolor{black}{we model the prior knowledge of the} code output \textcolor{black}{as} a Gaussian process $Z(\bm{x},t_u)$ with $\bm{x}\in Q$ and $u = 1, \hdots, N_t$ \textcolor{black}{with a} covariance function $C$ \textcolor{black}{given by} \cref{eq::CovFuncTenso} and \textcolor{black}{a} mean function $\mu$ \textcolor{black}{given by} \cref{eq::tredTensCov}.
We focus our attention to the case $N_t>N_x$.
We assume that the covariance structure can be decomposed into two different functions representing the correlation in $\bm{x}$ and the correlation in $t$.
If we choose well both functions, the kriging calculation is possible \cite{perrin2019adaptive,rougier2008efficient}.

\textcolor{black}{In the following we present a simplification of the method proposed in \cite{perrin2019adaptive}.
The} a priori $\R^{N_t}$-valued mean function \textcolor{black}{is assumed to be of the form}:
\begin{equation}
    \mu(\bm{x})=Bf(\bm{x})
    \label{eq::tredTensCov}
\end{equation}
where $f(\bm{x})$ is a given $\R^M$-valued function
and $B \in \mathcal{M}_{N_t\times M}(\R)$ is to be estimated.
We define \textcolor{black}{by} $F$ the $N_x \times M$ matrix  $[f^T(\bm{x}^{(i)})]_{i=1,\ldots,N_x}$.

The a priori covariance function $C(t_u,t_{u^\prime},\bm{x},\bm{x}^\prime)$ can be expressed with  the $N_t\times N_t$ matrix $R_t$ and the correlation function $C_x : Q\times Q \rightarrow [0,1]$ with $C_x(\bm{x},\bm{x})=1$:
\begin{equation}
    \label{eq::CovFuncTenso}
    C(t_u,t_{u^\prime},\bm{x},\bm{x}^\prime) = R_t(t_u,t_{u^\prime})C_x(\bm{x},\bm{x}^\prime) .
\end{equation}
\textcolor{black}{The covariance in time is expressed as a matrix because the temporal grid is finite and fixed.}
The covariance "matrix" (here a tensor) of $\left(Z(\bm{x}^{(j)},t_u)\right)_{\substack{u=1,\ldots,N_t\\ j=1,\ldots,N_x}}$ is 
\begin{equation}\label{matcorsep}
    R = R_t \otimes R_x,
\end{equation}
with $\left(R_x\right)_{k,l} =  C_x(\bm{x}^{(k)},\bm{x}^{(l)})$ $k,l=1,\dots,N_x$.

If $R_x$ and $R_t$ are not singular, then
 the a posteriori distribution of the $\R^{N_t}$-valued process $Z$ given the covariance functions and \textcolor{black}{the} observations \textcolor{black}{$\bm{Z}_\text{obs}$} is Gaussian:
\begin{equation}\label{processtimeserie}
 \left(Z(\bm{x},t_u)\right)_{u=1,\ldots,N_t} | R_t,C_x,\bm{Z}_\text{obs} \sim \mathcal{GP}(\mu_\star(\bm{x}) , R_\star(\bm{x},\bm{x}')R_t),
\end{equation}
with the $N_t$-dimensional posterior mean:
\begin{equation}
    \label{eq::meanTensCov}
    \mu_\star(\bm{x}) = \bm{Z}_\text{obs} R_x^{-1}r_x(\bm{x}) + B_\star u(\bm{x}) 
\end{equation}
where $r_x(\bm{x})$ is the $N_x$-dimensional vector $\left(C_x(\bm{x},\bm{x}^{(j)})\right)_{j=1,\ldots,N_x}$.
The posterior covariance function $R_\star (\bm{x},\bm{x}^\prime)$ is :
\begin{equation}
    \label{eq::varTensCov}
        R_\star(\bm{x},\bm{x}') = c_\star(\bm{x},\bm{x}')\left(1 +v_\star(\bm{x},\bm{x}')\right) .
\end{equation}
The functions that are used in the regression are
\begin{equation}
    \label{eq::TensCovpardev}
\left\{
\begin{array}{l}
u(\bm{x}) = f(\bm{x}) - F^TR_x^{-1}r_x(\bm{x}) \\
c_\star(\bm{x},\bm{x}^\prime) = C_x(\bm{x},\bm{x}^\prime) - r_x(\bm{x})^T R^{-1}_x r_x(\bm{x}')\\
v_\star(\bm{x},\bm{x}') = u(\bm{x})^T(F^TR_x^{-1}F)^{-1}u(\bm{x}')c_\star^{-1}(\bm{x},\bm{x}^\prime) \\
\end{array}
\right. ,
\end{equation}
and
\begin{equation}
    \label{eq::postEqAtemps}
    B_\star = \bm{Z}_\text{obs}R_x^{-1}F(F^TR_x^{-1}F)^{-1}.\\
\end{equation}

The correlation function $C_x$ is assumed to be a Matérn \textcolor{black}{$\frac{5}{2}$} kernel with a tensorized form\textcolor{black}{, see \cite[Chapter 4]{williams2006gaussian}}:
\begin{equation}
    C_x(\bm{x},\bm{x}^\prime)= \prod_{i=1}^{d}\left(1+\frac{\sqrt{5} | {x_i}- {x_i}^\prime |}{\ell_{x_i}}+ \frac{5| {x_i}-{x_i}^\prime |^2}{3\ell_{x_i}^2}\right) \exp\left(-\frac{\sqrt{5}| {x_i}-{x_i}^\prime|}{\ell_{x_i}}\right),
\end{equation}
with $\bm{\ell_x}=(\ell_{x_1},\ldots,\ell_{x_d})$ the vector of correlation lengths.
Other choices are of course possible.
$R_t$ is estimated using $R_x^{-1}$ and the observations $\bm{Z}_{\text{obs}}$ by maximum likelihood, as in \cite{perrin2019adaptive}:  
\begin{equation}
    \label{eq::CompRtcov}
    \widehat{R_t} = \frac{1}{N_x}\left(\bm{Z}_{\text{obs}} - \bm{\hat{Z}}\right) R_x^{-1} \left(\bm{Z}_{\text{obs}} - \bm{\hat{Z}}\right)^T,
\end{equation}
with $\bm{\hat{Z}}$ is the $N_t \times N_x$ matrix of empirical means $\hat{Z}_{u,i}=\frac{1}{N_x}\sum_{j=1}^{N_x}\left(\bm{Z}_\text{obs}\right)_{u,j}$, $\forall i=1, \ldots, N_x$ and $u=1,\ldots,N_t$.  

It remains only to estimate the vector of correlation lengths $\bm{\ell_x} = \left(\ell_{x_1},\dots,\ell_{x_d}\right)$ to determine the function $C_x$.
As presented in \cite{perrin2019adaptive}, the maximum likelihood estimation is not well defined for $\bm{\ell_x}$.
Indeed the approximation of $R_t$ by \cref{eq::CompRtcov}  is singular because $N_x<N_t$.
In fact we do not need to inver\textcolor{black}{t} $R_t$ as seen in \Cref{eq::postEqAtemps,processtimeserie,eq::meanTensCov,eq::varTensCov,eq::TensCovpardev}.
The method generally used to estimate the correlation lengths is cross-validation and, in our case, Leave-One-Out (LOO).
The LOO mean square error that needs to be minimized is:
\begin{equation}
    \label{eq::LOOerrnorm}
    \varepsilon^2(\bm{\ell_x}) = \sum_{k=1}^{N_x} \Vert \mu_\star^{(-k)}(\bm{x}^{(k)} |\bm{Z}_\text{obs}^{(-k)},\bm{l_x}) - \bm{Z}_\text{obs}(\bm{x}^{(k)})\Vert^2,
\end{equation}
where $\mu_\star^{(-k)}(\bm{x}^{(k)} |\bm{Z}_\text{obs}^{(-k)},\bm{l_x})$ is the $\R^{N_t}$-valued prediction mean obtained with the correlation length vector $\bm{l_x}$, using all observations except the $k$-th, at the point $\bm{x}^{(k)}$ \textcolor{black}{and $\Vert \cdot \Vert $ is the Euclidean norm in $\R^{Nt}$}.
We can use an expression of $\varepsilon^2(\bm{l_x})$ that does not require multiple regression, as in \cite{bachoc2013cross,dubrule1983cross}.
For more detail see \cref{asec::LOO}.

\section{AR(1) multi-fidelity model with projection}
\label{sec::AR1multiproj}
To carry out a \textcolor{black}{single}-fidelity regression for a code whose output is a time series we can use the method presented in \Cref{subsec::tensRegcov} or use output dimension reduction.
For dimension reduction, as in \cite{nanty2017uncertainty}, a basis is chosen.
The functional output is expanded onto this basis, the expansion is truncated and a surrogate model is built for each scalar-valued coefficient of the truncated expansion.
In our case, we deal with both multi-fidelity and time-series outputs\textcolor{black}{.}

One solution could be to use covariance tensorization in the multi-fidelity framework. This method requires the inversion of large covariance matrices.
However, the inversion methods \textcolor{black}{that are} possible in \textcolor{black}{a single-}fidelity \textcolor{black}{framework} become impossible \textcolor{black}{in a multi-fidelity framework} because the matrices are \textcolor{black}{then} too ill-conditioned to be inverted.

This leads us to introduce new methods.
The most naive method consists in starting again from the dimension reduction \textcolor{black}{technique} and to carry out the projection of the outputs of the two codes onto the same appropriate basis.
It is therefore possible to use the model \Cref{multikrigeq}.
The problem is that is not possible to define a basis that is optimal for both the high\textcolor{black}{-} and low\textcolor{black}{-}fidelity codes.
A basis estimated from the low-fidelity data is preferred as it is more robust thanks to the larger number of data.
Thus the loss of high-fidelity information is significant which leads us to introduce our new method \textcolor{black}{in} \Cref{sec::AR1multiTens}.
The originality of our approach is to keep the dimension reduction technique but also to use the residual high-fidelity data to carry out a \textcolor{black}{GP} regression with covariance tensorization.

\subsection{Model}
We recall that $\left(Z_\text{L}, Z_\text{H}\right)$ is a stochastic process.
The temporal grid is $\left\{ t_u \right\}_{u\in\left\{1,\hdots ,N_t\right\}}$.
$N_\text{F}$ observations are available for different values of $\bm{x}$ at the fidelity $\text{F}$.

Let $\bm{\Gamma}$ be an orthogonal $N_t\times N_t$ matrix.
The columns of $\bm{\Gamma}$, $\Gamma_i$, form an orthonormal basis of $\R^{N_t}$.
We assume that, given $\bm{\Gamma}$, \textcolor{black}{$\left( Z_\text{L}, Z_\text{H}\right)$ are Gaussian processes with covariance matrices that can be diagonalized on the basis formed by the columns of $\bm{\Gamma}$.}
\textcolor{black}{T}he processes $Z_\text{\textcolor{black}{F}}$ \textcolor{black}{can then be expanded as}:
\begin{align}\label{eq::decompBB}
Z_\text{L}(\bm{x},t_u) = & \sum_{i=1}^{N_t} A_{i,\text{L}}(\bm{x}) \Gamma_i(t_u),\\
\label{eq::decompBH}
Z_\text{H}(\bm{x},t_u) =& \sum_{i=1}^{N_t} A_{i,\text{H}}(\bm{x}) \Gamma_i(t_u),
\end{align}
where $\left(A_{i,\text{L}}(\bm{x}), A_{i,\text{H}}(\bm{x})\right)$ are Gaussian processes which are independent with respect to $i$, given $\bm{\Gamma}$.

Let $\left(\alpha_{i,\text{F}}(\bm{x})\right)_{i=1}^{N_t}$ be the $\R^{N_t}$-valued function:
\begin{equation}\label{eq::alphadefcomp}
\alpha_{i,\text{F}}(\bm{x}) = \sum_{u=1}^{N_t} z_{\text{F}}(\bm{x},t_u) \Gamma_i(t_u).
\end{equation}
We denote by $\alpha_{i}^{\text{F}}$ the $1\times N_F$ row vector $\left(\alpha_{i,\text{F}}(\bm{x}^{(j)})\right)_{j=1}^{N_F}$  that contains the available data.
The full data set is $\alpha=\left(\alpha^{\text{L}},\alpha^{\text{H}}\right)$.

Consequently we will use the method presented in \Cref{subsec::multiGPreg} given $\bm{\Gamma}$.
This leads us to the model presented in \Cref{eq::projModel}.
Given $\bm{\Gamma}$, $\forall i \in \left\{1, \hdots, N_t\right\}$,

\begin{equation}\label{eq::projModel}
    \left\{
    \begin{array}{rcl}
        A_{i,\text{H}}(\bm{x}) & = & \rho_{i,\text{L}}(\bm{x})\tilde{A}_{i,\text{L}}(\bm{x})+\delta_i(\bm{x})\\ 
        \tilde{A}_{i,\text{L}}(\bm{x}) & \perp & \delta_i(\bm{x})\\
        \rho_{i,\text{L}}(\bm{x}) & = & g^T_i(\bm{x})\beta_{\rho_\text{L},i}\\
    \end{array}
    \right. ,
\end{equation}
where:
$$ \left[ \delta_i(\bm{x}) | \bm{\Gamma}, \sigma_{i,H}, \beta_{i,H} \right] \sim \mathcal{GP}\left(f^T_{i,\text{H}}(\bm{x})\beta_{i,\text{H}},\sigma_{i,\text{H}}^2r_{i,\text{H}}(\bm{x},\bm{x}^\prime)\right),$$
and $ \tilde{A}_{i,\text{L}}(\bm{x})$ a Gaussian process conditioned by the values $\alpha^{\text{L}}$.
The distribution of $\tilde{A}_{i,\text{L}}(\bm{x})$ is the one of $\left[A_{i,\text{L}}(\bm{x})|\bm{\Gamma}, \mathcal{A}^{\text{L}}=\alpha^{\text{L}},\beta_{i,\text{L}},\sigma_{i,\text{L}}\right]$ where:

\begin{equation}
    \label{eq::lawAiLdispo}
    \left[ A_{i,\text{L}}(\bm{x}) | \bm{\Gamma}, \sigma_{i,L}, \beta_{i,L} \right] \sim \mathcal{GP}\left(f^T_{i,\text{L}}(\bm{x})\beta_{i,\text{L}},\sigma_{i,\text{L}}^2r_{i,\text{L}}(\bm{x},\bm{x}^\prime)\right).
\end{equation}

$g_i(\bm{x})$ \textcolor{black}{are} vectors of $q$ regression functions,
$f_{i,\text{F}}(\bm{x})$ are vectors of $p_\text{F}$ regression functions,
$r_{i,\text{F}}(\bm{x},\bm{x}^\prime)$ are correlation functions,
$\beta_ {i,\text{F}}$ \textcolor{black}{are} $p_\text{F}$-dimensional vectors,
$\beta_{\rho_\text{L},i}$ \textcolor{black}{are} $q$-dimensional vectors and $\sigma_{i,\text{F}}^2$ are \textcolor{black}{positive real numbers}.
For simplicity the regression functions $g_i$ and $f_{i,\text{F}}$ do not depend on $i$.

The model depends on $\bm{\Gamma}$, which is why we discuss in \Cref{ssec::BaseGamma} the choice of the basis.

\subsection{Basis}
\label{ssec::BaseGamma}
In this section we present different models for the random orthogonal matrix $\bm{\Gamma}$.
Its law depends on \textcolor{black}{the available information}.
If we have access to a lot of information based on the output of our code, we can use a Dirac distribution concentrated on one orthogonal matrix $\bm{\gamma}$ (it is a form of plug-in method).
In contrast, the least informative law is the Uniform Law, i.e. the Haar measure over the group of orthogonal matrices.
In order to make the best use of the available information, i.e. the known results of the code, an empirical law can be used. 

\subsubsection{Dirac distribution} 
\label{sssec::DiracLaw}

We can choose the distribution of the random matrix $\bm{\Gamma}$ as a Dirac distribution concentrated on a well chosen orthogonal matrix $\bm{\gamma}$.
This matrix is chosen when the basis is known \textcolor{black}{a priori} or if the basis can be efficiently estimated from the observed code outputs.
Motivated by \textcolor{black}{the} remark below \Cref{eq::decompBH}, the matrix $\bm{\gamma}$ can be computed using the singular value decomposition (SVD) of the code outputs.

The general idea is to choose subsets $\tilde{D}_\text{F} \subset D_\text{F} $ of size $\tilde{N}_\text{F}$ and to apply a SVD on the $N_t \times \left(\tilde{N}_\text{H}+\tilde{N}_\text{L}\right)$ matrix $\tilde{\bm{Z}}_\text{obs}$ that contains the observed values $\left(z_\text{H}(\bm{x},t_u)\right)_{\substack{u = 1,\ldots,N_t\\ \bm{x}\in \tilde{D}_\text{H}}}$ and $\left(z_\text{L}(\bm{x},t_u)\right)_{\substack{u = 1,\ldots,N_t\\ \bm{x}\in \tilde{D}_\text{L}}}$.
The SVD gives:
\begin{equation}
    \tilde{\bm{Z}}_\text{obs} = \bm{\tilde{U}} \bm{\tilde{\Lambda}} \bm{\tilde{V}}^T.
\end{equation}
The choice of $\bm{\gamma}$ is $\bm{\tilde{U}}$.

The first idea is to mix all available data, high- and low-fidelity: $\tilde{D}_\text{H}= D_\text{H}$ and $\tilde{D}_\text{L}=D_\text{L}$.
However, we \textcolor{black}{typically have} that $N_\text{L}\gg N_\text{H}$, so the basis is mainly built \textcolor{black}{from} the low-fidelity data.
In addition, the small \textcolor{black}{differences} in the data between high\textcolor{black}{-} and low\textcolor{black}{-}fidelity \textcolor{black}{code outputs} that would be useful to build the basis have negligible impact because they are overwhelmed by the low\textcolor{black}{-}fidelity data.
This method \textcolor{black}{is not appropriate in our framework}.

We have to choose between high\textcolor{black}{-} and low\textcolor{black}{-} fidelity.
High\textcolor{black}{-}fidelity has the advantage of being closer to the desired result.
However, it is also almost impossible to validate the chosen $\bm{\gamma}$ because the high-fidelity data size $N_\text{H}$  is small.
The low-fidelity data set is larger, hence the estimation of $\bm{\gamma}$ is more robust.
In order to choose $\bm{\gamma}$, we therefore suggest to use the low-fidelity data and to calculate the SVD with $\tilde{D}_\text{H}= \emptyset$ and $\tilde{D}_\text{L}=D_\text{L}$.

\subsubsection{Uniform distribution}

We can choose a random matrix \textcolor{black}{$\bm{\Gamma}$} using the Haar measure on the orthogonal group $O_{N_t}$, the group of $N_t \times N_t$ orthogonal matrices.
This is the Uniform Orthogonal Matrix Law.

To generate a random matrix from the Haar measure over $O_{N_t}$, one can first generate a $N_t \times N_t$ matrix with independent and identically distributed coefficients with the reduced normal distribution, then, apply the Gram-Schmidt process onto the matrix.
As shown in \cite{diaconis2005random}, this generator produces a random orthogonal matrix with the uniform orthogonal matrix law.
This method completely ignores the available data and is not appropriate in our framework.

\subsubsection{\textcolor{black}{Cross-Validation Based} distribution}
\label{sssec::ELprojbase}

The downside of the Dirac distribution is that the uncertainties on the basis \textcolor{black}{estimation} are not taken into account.
A \textcolor{black}{cross-validation based (CVB)} method to assess the uncertainty estimation is therefore considered.

The proposed method uses only the low-fidelity data because it is assumed that there are too few high-fidelity data to implement this method, so $\tilde{D}_\text{H}=\emptyset$.
For the construction of the basis we try to have different sets to evaluate the basis in order to have empirical estimates of the moments of the basis vectors.
Let $k$ be a fixed integer in $\{1, \hdots, N_\text{L}\}$.
Let $I= \{J_1, \hdots, J_k\}$ be a random set of $k$ elements in $\{1, \hdots, N_\text{L}\}$, with uniform distribution over the subsets of $k$ elements in $\{1,\hdots, N_\text{L}\}$.
The empirical distribution \textcolor{black}{for the matrix $\bm{\Gamma}$} is defined as follows: for any \textcolor{black}{bounded function} $f: O_{N_t} \rightarrow \R$,
\begin{equation}
    \label{eq::emplowGamma}
    \E\left[f(\bm{\Gamma}_I)\right]= \frac{1}{\binom{N_\text{L}}{k}}\sum_{\{j_1,\ldots,j_k\}\subset \{1,\ldots,N_\text{L}\}}f(\bm{\tilde{U}}_{[J_1,\ldots,J_k]})  ,
\end{equation}
where $\bm{\tilde{U}}_{[J_1,\ldots,J_k]}$ is the matrix of the left  singular vectors of the SVD of\\ $\left(z_\text{L}(\bm{x}^{(i)},t_u)\right)_{\substack{u\in\{1,\ldots,N_t\}\\ i\in \{1,\ldots,N_\text{L}\}\backslash\{J_1,\ldots,J_k\}}}$.
This distribution depends on the choice of $k$, that will  be discussed in \Cref{sec::Illustration}.

\subsection{Predictive mean and variance}
The goal of this section is to calculate the posterior distribution of $Z_\text{H}(\bm{x},t_u)$.
The problem can be split into two parts: the multi-fidelity regression of the basis coefficients knowing $\bm{\Gamma}$ and the integration with respect to the distribution of $\bm{\Gamma}$.
The Dirac and \textcolor{black}{CVB} distributions \textcolor{black}{described in \cref{ssec::BaseGamma}} can be used to define the law of $\bm{\Gamma}$.

\paragraph{Multi-fidelity surrogate modeling of the coefficients}
\label{sssec::MuFIRegProj}
By applying the model proposed in \Cref{subsec::multiGPreg} we can therefore deduce the prediction mean and variance.
\textcolor{black}{Their expressions are given in \cref{assec::MultiUnco}}

\subsubsection{Dirac law of \texorpdfstring{$\bm{\Gamma}$}{Gamma}}
Here we assume that the law of $\bm{\Gamma}$ is Dirac at $\bm{\gamma}$.
Consequently, the posterior distribution of $Z_\text{H}(\bm{x},t)$ is Gaussian.
In order to characterize the law of $Z_\text{H}(\bm{x},t_u)$ it is necessary and sufficient to compute its mean and variance.

\paragraph{Mean:}
The posterior mean is:
\begin{equation}
    \label{eq::MeanZparaDirac}
    \E\left[Z_\text{H}(\bm{x},t_u) | \mathcal{A}=\alpha\right] =
    \sum_{i=1}^{N_t} \gamma_i(t_u) \E\left[A_{i,\text{H}}(\bm{x}) | \mathcal{A}=\alpha\right],
\end{equation}
where the expectation $\E\left[A_{i,\text{H}}(\bm{x}) | \mathcal{A}=\alpha\right]$ is given by \cref{eq::meanautomodeleq}.

\paragraph{Variance:}
The posterior variance:
\begin{equation}
    \label{eq::VarZparaDirac}
    \V\left[ Z_{\text{H}}(\bm{x},t_u) | \mathcal{A}=\alpha\right] = 
    \sum_{i=1}^{N_t}\gamma_i^2(t_u) \V\left[A_{i,\text{H}}(\bm{x})|\mathcal{A}=\alpha\right] ,
\end{equation}
where the variance $\V\left[A_{i,\text{H}}(\bm{x})|\mathcal{A}=\alpha\right]$ is given by \cref{eq::varautomodeleq}. 

\subsubsection{\textcolor{black}{CVB} law of \texorpdfstring{$\boldsymbol{\Gamma}$}{Gamma}}
\label{sssec::projArbGamma}
Because the law is different from Dirac the posterior distribution of $Z_\text{H}(\bm{x},t)$ is not Gaussian anymore.
However, we can characterize the posterior mean and the variance of $Z_\text{H}(\bm{x},t_u)$.

We denote $\E_\alpha\left[\cdot\right]= \E\left[\cdot|\mathcal{A}=\alpha\right]$, $\V_\alpha\left[\cdot\right]= \V\left[\cdot|\mathcal{A}=\alpha\right]$,
$\E_{\bm{Z}_\text{obs}}\left[\cdot\right]= \E\left[\cdot|\mathcal{Z}=\bm{Z}_\text{obs}\right]$ and $\V_{\bm{Z}_\text{obs}}\left[\cdot\right]= \V\left[\cdot|\mathcal{Z}=\bm{Z}_\text{obs}\right]$.

\paragraph{Mean}
The linearity of the expectation and the law of total expectation give:
\begin{equation}
    \label{eq::MeanZpara}
    \E_\alpha\left[Z_\text{H}(\bm{x},t_u)\right] = \sum_{i=1}^{N_t}\E_\alpha{\left[\Gamma_i(t_u)\E_\alpha{\left[A_{i,\text{H}}(\bm{x}) |\bm{\Gamma} \right]}\right]},
\end{equation}
where the expectation $\E_\alpha{\left[A_{i,\text{H}}(\bm{x}) |\bm{\Gamma} \right]}$ is given by \Cref{eq::meanautomodeleq}.

\paragraph{Variance}
The law of total variance gives :
\begin{equation}
    \label{eq::VarZpara}
        \V_\alpha\left[ Z_{\text{H}}(\bm{x},t_u) \right] =
        \V_\alpha\left[\E_\alpha\left[ Z_\text{H}(\bm{x},t_u)|\bm{\Gamma}\right]\right]
        + \E_\alpha\left[\V_\alpha\left[ Z_\text{H}(\bm{x},t_u)|\bm{\Gamma}\right]\right]
\end{equation}
\textcolor{black}{By \cref{assec:varproj} we get:}
\begin{equation}
    \label{eq::VarZparaFull}
    \begin{array}{c}
        \V_\alpha\left[ Z_{\text{H}}(\bm{x},t_u) \right] = 
        \sum_{i=1}^{N_t} \V_\alpha\left[\Gamma_i (t_u)\E_\alpha \left[A_{i,\text{H}}(\bm{x})|\bm{\Gamma}\right]\right] \\
        + \sum_{i,j=1,i\neq j}^{N_t} \Cov_\alpha (\Gamma_i (t_u)\E_\alpha\left[A_{i,\text{H}}(\bm{x})|\bm{\Gamma}\right], \Gamma_j(t_u) \E_\alpha\left[A_{j,\text{H}}(\bm{x})|\bm{\Gamma} \right])\\
        + \sum_{i=1}^{N_t}\E_\alpha\left[\Gamma_i^2(t_u) \V_\alpha\left[A_{i,\text{H}}(\bm{x})|\bm{\Gamma}\right]\right] \\
    \end{array}.
\end{equation}

\Cref{eq::VarZparaFull,eq::MeanZpara} are combinations of expectations of explicit functions of $\bm{\Gamma}$.
We can compute the result using our knowledge on the law of $\bm{\Gamma}$.
The expectation of a function of $\bm{\Gamma}$ is given by \Cref{eq::emplowGamma}.

\subsection{Truncation}
\label{sssec::truncProjec}

There is a problem with \textcolor{black}{the} surrogate modeling of the coefficients of the decomposition with indices larger than $N_\text{L}$.
Indeed, we typically have $N_\text{L} < N_t$ so the vectors $\Gamma_i$ with indices larger than $N_\text{L}$ of the basis are randomly constructed, which is not suitable for building surrogate models.
To solve this problem, it is possible to truncate the sum.
Only the first $N$ coefficients, with  $N\leq N_\text{L}$ are calculated.
\textcolor{black}{This would be reasonable if} the contributions of the terms $A_{i,\text{H}}(\bm{x})\Gamma_i(t_u)$ for $i>N$ \textcolor{black}{were} negligible.
However it turns out that these terms are often not collectively negligible (see \Cref{sec::Illustration}) and the truncation method does not achieve a good bias-variance trade-off even when optimizing with respect to $N$ (by a cross validation procedure for instance).
The high- and low-fidelity outputs do not necessarily have the same forms.
Thus it is possible that an important part of the high-fidelity code is neglected because it is not taken into account by the sub-space spanned by $\{\Gamma_i\}_{i\leq N}$.
We will therefore propose in the next section an original method to tackle this problem.

\section{AR(1) multi-fidelity model with tensorized covariance and \textcolor{black}{p}rojection}
\label{sec::AR1multiTens}
The naive method presented in \Cref{sec::AR1multiproj} has many flaws.
It leaves part of the output untreated by regression and the variance is underestimated.
The major problem with the solution we propose in \Cref{sec::AR1multiproj}, is that we typically have $N_\text{L}<N_t$.
Consequently for $i>N_\text{L}$ \textcolor{black}{the} vectors $\Gamma_i$ of the basis do not represent the typical variation\textcolor{black}{s} of $Z_\text{H}$.
We should find an appropriate way to predict the law of the projection of $Z_\text{H}$ on $\operatorname{span}\left\{\Gamma_i, i>N_\text{L}\right\}$.

One interesting approach is to apply the covariance tensorization method to the orthogonal part.
This idea is to address the last terms of the expression collectively through a \textcolor{black}{GP} model with tensorized covariance structure.
This allows us  to split the problem into two parts. 
The first part is to compute (as presented in the previous section) the first $N$ terms of the expansion of $Z_\text{H}$ onto the basis by using a co-kriging approach.
The second part is to compute the projection of $Z_\text{H}$ onto the orthogonal space by using a kriging approach with tensorized covariance.
The choice of the optimal $N$ will be carried out by a K-fold cross-validation method.

\subsection{Decomposition}
The proposed method is based on the decomposition of the outputs, as presented in \Cref{eq::decompBB}.

\paragraph{Projection}
Let $N$ be an integer, smaller than the time dimension $N_t$.
Let $\bm{\Gamma}$ be an orthogonal matrix, as presented in \Cref{ssec::BaseGamma}, the columns of $\bm{\Gamma}$ are $\{\Gamma_i\}_{i=1,\cdots,N_t}$.
As discussed in \Cref{sssec::truncProjec} the full computation of all $N_t$ surrogate models may not give good results.
This leads to the idea of the introduction of the orthogonal subspaces $S_N^\parallel$ and $S_N^\perp$, where $S_N^\parallel=\operatorname{span}\left\{\Gamma_1,\ldots,\Gamma_N\right\}$ and $S_N^\perp=\operatorname{span}\left\{\Gamma_{N+1},\ldots,\Gamma_{N_t}\right\}$.

With a given basis $\bm{\Gamma}$ it is possible to decompose the code outputs.
The decomposition over the basis $\bm{\Gamma}$ gives us coefficients.
We decompose $Z_{\text{H}}$ and $Z_{\text{L}}$ over the subspace $S_N^\parallel$.
The rests are denoted $Z^\perp_{\text{H}}$ and $Z^\perp_{\text{L}}$.
Consequently, we get:

\begin{equation}\label{eq::SVDlowFi}
Z_{\text{L}}(\bm{x},t_u) = Z_{\text{L}}^{\parallel}(\bm{x},t_u) + Z_{\text{L}}^{\perp}(\bm{x},t_u)
= \sum_{i=1}^N A_{i,\text{L}}(\bm{x})\Gamma_i(t_u) + Z^\bot_\text{L}(\bm{x},t_u)
\end{equation}and
\begin{equation}\label{eq::SVDhighFi}
Z_{\text{H}}(\bm{x},t_u) = Z_{\text{H}}^{\parallel}(\bm{x},t_u) + Z_{\text{H}}^{\perp}(\bm{x},t_u)
= \sum_{i=1}^N A_{i,\text{H}}(\bm{x})\Gamma_i(t_u) + Z_{\text{H}}^{\perp}(\bm{x},t_u).
\end{equation}

We are able to describe the code outputs with the basis $\bm{\Gamma}_N=\left\{\Gamma_i\right\}_{i=1,\ldots,N}$ of $S_N$, the coefficients $A_{i,\text{H}}$ and $A_{i,\text{L}}$, and the orthogonal parts $Z^\bot_\text{L}$ and $Z^\bot_\text{H}$.
We denote by $\alpha_{i,\text{H}}$ and $\alpha_{i,\text{L}}$  the available data sets, the full set is called $\alpha$.
The expression of $\alpha_{i,{\text{F}}}$ is given by \Cref{eq::alphadefcomp}.

We will use the method presented in \Cref{subsec::multiGPreg} for all $i\leq N$.
Given $\bm{\Gamma}$,

\begin{equation}\label{eq::projModelmulti}
    \left\{
    \begin{array}{rcl}
        A_{i,\text{H}}(\bm{x}) & = & \rho_{i,\text{L}}({\bm{x}})\tilde{A}_{i,\text{L}}(\bm{x})+\delta_i(\bm{x})\\ 
        \tilde{A}_{i,\text{L}}(\bm{x}) & \perp & \delta_i(\bm{x})\\
        \rho_{i,\text{L}}(\bm{x}) & = & g^T_i(\bm{x})\beta_{i,\rho_\text{L}}\\
    \end{array}
    \right.,
\end{equation}
where:
$$ \left[ \delta_i(\bm{x})| \bm{\Gamma}, \beta_{i,H}, \sigma_{i,H} \right] \sim \mathcal{GP}\left(f^T_{i,\text{H}}(\bm{x})\beta_{i,\text{H}},\sigma_{i,\text{H}}^2r_{i,\text{H}}(\bm{x},\bm{x}^\prime)\right),$$
and $\tilde{A}_{i,\text{L}}(\bm{x})$ is a Gaussian process conditioned by $\alpha^{\text{L}}$.
Its distribution is the one of $\left[A_{i,\text{L}}(\bm{x})|\bm{\Gamma}, \mathcal{A}^{\text{L}}=\alpha^{\text{L}},\beta_\text{L},\sigma_\text{L}\right]$ where
\textcolor{black}{the law of $\left[ A_{i,\text{L}}(\bm{x}) | \bm{\Gamma}, \beta_{i,L}, \sigma_{i,L} \right]$ is of the form \cref{eq::lawAiLdispo}.}

$$ \left[ A_{i,\text{L}}(\bm{x}) | \bm{\Gamma}, \beta_{i,L}, \sigma_{i,L} \right] \sim \mathcal{GP}\left(f^T_{i,\text{L}}(\bm{x})\beta_{i,\text{L}},\sigma_{i,\text{L}}^2r_{i,\text{L}}(\bm{x},\bm{x}^\prime)\right).$$

$g_i$ are vectors of $q$ regression functions,
$f_{i,\text{F}}(\bm{x})$ are vectors of $p_\text{F}$ regression functions,
$r_{i,\text{F}}(\bm{x},\bm{x}^\prime)$ are correlation functions,
$\beta_ {i,\text{F}}$ are $p_\text{F}$-dimensional vectors,
$\beta_{i,\rho_\text{L}}$ are $q$-dimensional vectors
and $\sigma_{i,\text{F}}^2$ are positive real numbers.

For the orthogonal part projected onto $S_N^\perp$ the method is different.
The hypothesis is that the projection $Z^\perp_\text{L}(\bm{x},t_u)$ of $Z_\text{L}(\bm{x},t_u)$ has a negligible influence on the projection $Z^\perp_\text{H}(\bm{x},t_u)$ of $Z_\text{H}(\bm{x},t_u)$.
Our assumption is that $Z_\text{H}^\perp(\bm{x},t_u)$ is a Gaussian process with a tensorized covariance.
The method we will use on $Z_\text{H}^\perp(\bm{x},t_u)$ is described in \Cref{subsec::tensRegcov}.

\textcolor{black}{Note that t}he value $N=0$ corresponds to full \textcolor{black}{single}-fidelity, in this case we use only \textcolor{black}{GP} regression with covariance tensorization as in \cref{subsec::tensRegcov}.
For $N=N_\text{L}$ the dimension reduction is minimal and co-kriging is applied to all pairs $\left(A_{i,\text{L}},A_{i,\text{H}}\right)$ for $i\leq N_\text{L}$.
We will see \textcolor{black}{in \cref{sec::Illustration}} that the optimal $N$ is in fact positive but smaller than $N_\text{L}$.

\subsection{Predictive mean and covariance}
\label{subsec::CompPredTensor}
In this section we first make a quick reminder of the methods presented in \Cref{sec::Regressiontools}.
Moreover, with different assumptions \textcolor{black}{about} the law of $\bm{\Gamma}$ we present the regression using the model and the data.

\paragraph{Multi-fidelity of coefficients}
As in \cref{sssec::MuFIRegProj} we compute the $N$ multi-fidelity models of the first $N$ coefficients of the expansion of the code output given $\bm{\Gamma}$.
If we apply the method proposed in \cref{subsec::multiGPreg} we can therefore deduce the prediction mean and variance\textcolor{black}{, as in \cref{sssec::MuFIRegProj}}.

\paragraph{Tensorized covariance regression}
The orthogonal part of the regression is computed using the method presented in \cref{subsec::tensRegcov}.
The adaptation is that the regression must be \textcolor{black}{carried out} in subspace $S_N^\perp$ given $\bm{\Gamma}$,
\begin{equation}
    \label{eq::alphaortho}
    \bm{Z}_\text{obs}^\perp = \bm{Z}_\text{obs}
    - \left( \sum_{i=1}^{N} \alpha_{i,\text{H}}(\bm{x})\Gamma_i(t_u)\right)_{\substack{u=1,\cdots,N_t\\\bm{x} \in D_\text{H}}},
\end{equation}
where $\alpha_{i,\text{H}}$ is given by \Cref{eq::alphadefcomp}.

This does not have any consequence on the $\bm{x}$ part but only on the $t$ part.
Contrarily to $Z_\text{H}^\parallel(\bm{x},t_u)$ only one surrogate model is needed for $Z_\text{H}^\perp(\bm{x},t_u)$.
\textcolor{black}{The detail of how we can deal with $Z_\text{H}^\perp(\bm{x},t_u)$ and $Z_\text{H}^\parallel(\bm{x},t_u)$ is explained in \cref{asec::TensOrt}.}

\subsubsection{Dirac law of \texorpdfstring{$\bm{\Gamma}$}{Gamma}}
Here we assume that $\bm{\Gamma}$ is known and its distribution is Dirac at $\bm{\gamma}$.
Consequently, as in \Cref{sec::AR1multiproj}, $Z_\text{H}(\bm{x},t)$ is a Gaussian process by linear combination of independent Gaussian processes.
Its posterior distribution is completely determined if we can evaluate its mean and covariance.

\subparagraph{Mean}
The $\Gamma_i(t_u)$\textcolor{black}{' s are} constant and equal to $\gamma_i(t_u)$. Consequently: 
\begin{equation}
    \E_\alpha{\left[Z_{\text{H}}^{\parallel}(\bm{x},t_u)|N\right]} = \sum_{i=1}^N \E_\alpha{\left[A_{i,\text{H}}(\bm{x})\right]} \gamma_i(t_u),
\end{equation}
and
\begin{equation}
    \label{eq::DiracMeanTensor}
    \begin{array}{c}
        \E_{\bm{Z}_\text{obs}}{\left[Z_{\text{H}}(\bm{x},t_u)|N,\bm{\ell_x}\right]}
        = \sum_{i=1}^N \E_\alpha{\left[A_{i,\text{H}}(\bm{x})\right]} \Gamma_i(t_u)
        + \E_{\bm{Z}_\text{obs}}^\perp{\left[Z_{\text{H}}^{\perp}(\bm{x},t_u) | N,\bm{\ell_x}\right]}
    \end{array},
\end{equation}
where $\E_{\alpha}\left[A_{i,\text{H}}(\bm{x})\right]$ is given by (\ref{eq::meanautomodeleq}) and $\E_{\bm{Z}_\text{obs}}^\perp{\left[Z_{\text{H}}^{\perp}(\bm{x},t_u) | N,\bm{\ell_x}\right]}$ by 
(\ref{eq::meanTensCov}).

\subparagraph{Variance} 
The formula of the variance is:
\begin{equation}
    \V_\alpha\left[A_{i,\text{H}}(\bm{x})\Gamma_i(t) \right] =  \V_\alpha \left[A_{i,\text{H}}(\bm{x})\right]\gamma_i(t)^2 .
    \label{eq::varDiraclawtens}
\end{equation}

The uncorrelation of the coefficients $A_{i,\text{F}}(\bm{x})$ gives $\operatorname{Cov}_\alpha\left[A_{i,\text{H}}(\bm{x}),A_{j,\text{H}}(\bm{x})\right]=0$, for $i\neq j$ and $\operatorname{Cov}_\alpha\left[A_{i,\text{H}}(\bm{x})\gamma_i(t_u),Z^\bot_\text{H}(\bm{x},t_u)\right]=0$.
The expression of the variance becomes simple:

 \begin{equation}\label{eq::fullvarsimple}
\begin{array}{c}
 \V_{\bm{Z}_\text{obs}}\left[ Z_{\text{H}}(\bm{x},t_u) |N,\bm{\ell_x}\right] = 
 \sum_{i=1}^N \V_{\alpha}\left[A_{i,\text{H}}(\bm{x})\right]\gamma_i(t_u)^2 
 + \V_{\bm{Z}_\text{obs}}^\perp\left[Z^\bot_\text{H}(\bm{x},t_u)|N,\bm{\ell_x}\right], \end{array}
\end{equation}
where $\V_\alpha\left[A_{i,\text{H}}(\bm{x})\right]$ is given by (\ref{eq::varautomodeleq}) and
 $\V_{\bm{Z}_\text{obs}}^\perp\left[Z^\bot_\text{H}(\bm{x},t_u) | N, \bm{\ell_x} \right]$ is given by (\ref{processtimeserie}). 

\subsubsection{\textcolor{red}{CVB} law of \texorpdfstring{$\bm{\Gamma}$}{Gamma}}

\textcolor{black}{The posterior distribution of} $Z_\text{H}(\bm{x},t)$ is not Gaussian \textcolor{black}{anymore}.
However to predict the output of the high-fidelity code and to quantify the prediction uncertainty we are able to compute the posterior mean and variance of $Z_\text{H}(\bm{x},t)$.

\subparagraph{Mean}
We can decompose the process into two parts:
\begin{equation}
 Z_{\text{H}}(\bm{x},t_u) =  Z_{\text{H}}^{\parallel}(\bm{x},t_u) + Z_{\text{H}}^{\perp}(\bm{x},t_u).
\end{equation}

The linearity of the expectation gives us:
%
\begin{equation}
    \E_{\bm{Z}_\text{obs}}{\left[Z_{\text{H}}(\bm{x},t_u)|N,\bm{\ell_x}\right]} 
    =\sum_{i=1}^N \E_{\bm{Z}_\text{obs}}{\left[A_{i,\text{H}}(\bm{x}) \Gamma_i(t_u)\right]}
    + \E_{\bm{Z}_\text{obs}}^\perp\left[Z_{\text{H}}^{\perp}(\bm{x},t)|N,\bm{\ell_x}\right] .
\end{equation}
The theorem of total expectation gives us:
\begin{equation}
\label{eq::fullorthoesp}
\E_{\bm{Z}_\text{obs}}{\left[Z_{\text{H}}^{\parallel}(\bm{x},t_u)|N\right]} = \sum_{i=1}^N \E_{\bm{Z}_\text{obs}}{\left[\Gamma_i(t_u)\E_{\alpha}\left[A_{i,\text{H}}(\bm{x}) |\bm{\Gamma}\right]\right]},
\end{equation}
and therefore,
\begin{align}
\nonumber
        \E_{\bm{Z}_\text{obs}}\left[Z_\text{H}(\bm{x},t_u)|N,\bm{\ell_x}\right]
        =& \sum_{i=1}^N \E_{\bm{Z}_\text{obs}}{\left[\Gamma_i(t_u)\E_\alpha{\left[A_{i,\text{H}}(\bm{x}) |\bm{\Gamma} \right]}\right]}\\
       & + \E_{\bm{Z}_\text{obs}}{\left[\E_{\bm{Z}_\text{obs}}^\perp{\left[Z_{\text{H}}^{\perp}(\bm{x},t_u)|N,\bm{\ell_x},\bm{\Gamma}\right]}|N,\bm{\ell_x}\right]} .
    \label{eq::meanArbTens}
\end{align}
where $ \E_\alpha{\left[A_{i,\text{H}}(\bm{x}) |\bm{\Gamma}\right]} $ is given by \Cref{eq::meanautomodeleq} and $ \V_\alpha\left[Z^\perp_{\text{H}}(\bm{x})|\bm{\Gamma}\right] $ is given by \Cref{eq::meanTensCov}.
\Cref{eq::meanArbTens} is a combination of expectations of explicit functions of $\bm{\Gamma}$, which can be computed by \Cref{eq::emplowGamma}.

\subparagraph{Variance}

The theorem of the total variance gives us:
\begin{equation}
    \begin{array}{c}
        \V_{\bm{Z}_\text{obs}}\left[Z_{\text{H}}\left(\bm{x},t_u\right)|N,\bm{\ell_x}\right] =
        \V_{\bm{Z}_\text{obs}}\left[\E_{\bm{Z}_\text{obs}}\left[Z_{\text{H}}\left(\bm{x},t_u\right)|\bm{\Gamma},N,\bm{\ell_x}\right]|N,\bm{\ell_x}\right]\\
        + \E_{\bm{Z}_\text{obs}}\left[\V_{\bm{Z}_\text{obs}}\left[Z_{\text{H}}\left(\bm{x},t_u\right)|\bm{\Gamma},N,\bm{\ell_x}\right]|N,\bm{\ell_x}\right]
    \end{array}.
    \label{eq::NCalInterVar}
\end{equation}
\textcolor{black}{By \cref{assec:vartenso} we get:}
\begin{equation}
    \label{eq::VarFullortho}
    \begin{array}{c}
        \V_{\bm{Z}_\text{obs}}\left[Z_{\text{H}}\left(\bm{x},t_u\right)|N,\bm{\ell_x}\right] =
        \V_{\bm{Z}_\text{obs}}\left[\E_{\bm{Z}_\text{obs}}^\perp\left[ Z_\text{H}^\perp(\bm{x},t_u)|\bm{\Gamma},N,\bm{\ell_x}\right]|N,\bm{\ell_x}\right]\\
        + \E_{\bm{Z}_\text{obs}}\left[\V_{\alpha}\left[Z_\text{H}^\perp\left(\bm{x},t_u\right)|\bm{\Gamma}, N, \bm{\ell_x}\right]| N, \bm{\ell_x}\right]\\
        + \sum_{i=1}^N \V_{\bm{Z}_\text{obs}}\left[\Gamma_i(t_u)\E_{\alpha}\left[A_{i,\text{H}}(\bm{x})|\bm{\Gamma}\right]\right]\\ 
        + \sum_{i=1}^{N}\E_{\bm{Z}_\text{obs}}\left[\Gamma_i(t_u)^2\V_\alpha\left[A_i\left(\bm{x}\right)|\bm{\Gamma}\right]\right]\\
        + \sum_{i,j=1; i\neq j}^N \text{Cov}_{\bm{Z}_\text{obs}}\left[\Gamma_i(t_u)\E_\alpha\left[A_{i,\text{H}}(\bm{x})|\bm{\Gamma}\right],\Gamma_j(t_u)\E_\alpha\left[A_{j,\text{H}}(\bm{x})|\bm{\Gamma}\right]\right]\\
        + 2\sum_{i=1}^N \text{Cov}_{\bm{Z}_\text{obs}}\left[\Gamma_i(t_u)\E_\alpha\left[A_{i,\text{H}}(\bm{x})|\bm{\Gamma}\right], \E_{\bm{Z}_\text{obs}}^\perp\left[Z_\text{H}^\perp(\bm{x},t_u)|\bm{\Gamma}, N, \bm{\ell_x}\right]| N,  \bm{\ell_x}\right]\\
    \end{array}.
\end{equation}
where $ \V_\alpha{\left[A_{i,\text{H}}(\bm{x}) |\bm{\Gamma}\right]} $ is given by \Cref{eq::varautomodeleq} and $ \V_\alpha\left[Z_\text{H}^\perp(\bm{x})|\bm{\Gamma}\right] $ is given by \Cref{eq::varTensCov}.
\Cref{eq::VarFullortho} is a combination of expectations and variances of explicit functions of $\bm{\Gamma}$, which can be computed by \Cref{eq::emplowGamma}.

\subsection{Effective dimension}
For the formula\textcolor{black}{s} in \Cref{subsec::CompPredTensor} to be valid, $N$ must be fixed.
We may choose $N$ by a knowledge on the physical system or on the code but it is impossible in most cases due to the high/low-fidelity differences.
The best solution is generally to determine $N$ by a K-fold cross validation procedure.

The criterium that we choose to maximize is:
\begin{equation}
    Q^2_{N}(t_u) = 1 - 
    \frac{\sum_{k=1}^{N_\text{H}}\left(z_\text{H}(\bm{x}^{(k)},t_u)-\E\left[Z_\text{H}(\bm{x}^{(k)},t_u)|\bm{\Gamma}, N, \bm{\ell_x}, \bm{Z}_\text{obs}^{(-k)} \right]\right)^2}{N_\text{H}\V\left[z_\text{H}(D_\text{H},t_u)\right]},
\end{equation}
where $\V\left[z_\text{H}(D_\text{H},t_u)\right]$ is the empirical variance of the observed values:
$$
\V\left[z_\text{H}(D_\text{H},t_u)\right] 
=
\frac{1}{N_\text{H}} \sum_{k=1}^{N_\text{H}} z_\text{H}(\bm{x}^{(k)},t_u)^2  - 
\Big( \frac{1}{N_\text{H}} \sum_{k=1}^{N_\text{H}} z_\text{H}(\bm{x}^{(k)},t_u) \Big)^2 .
$$
The procedure we propose starts with the dimension $0$.
For the case $N=0$ the surrogate model depends only on high\textcolor{black}{-}fidelity regression.

\begin{itemize}
    \item[-] We compute the surrogate model for all $N=0, \ldots, N_\text{L}$
    \item[-] We calculate the mean in $t_u$ of $Q^2_{N}(t_u)$:
    $$ \hat{Q}_N^2 = \frac{1}{N_t}\sum_{u=1}^{N_t}Q_{N}^{2}(t_u)$$
\end{itemize}
\textcolor{black}{W}e compare the $\hat{Q}^2_{N}$ values and the value $N$ with the largest $\hat{Q}^2_{N}$ is chosen.
\textcolor{black}{In order to evaluate the surrogate model in the next section, we compute $\hat{Q}^2 = \max_N \hat{Q}^2_N$.}

%
%

\section{Illustration: double pendulum simulator}
\label{sec::Illustration}
The purpose of this section is to apply the methods proposed in the previous sections to a mechanical  example.
The example is based on a simulator of a pendulum attached to a spring-mass system.
We have two codes: the high-fidelity code numerically solves Newton's equation.
The low-fidelity code simplifies the equation, by linearisation for small angles of the pendulum motion, and solves the system.

\subsection{Characteristics of the outputs}
\begin{figure}
    \centering
    \includegraphics[height=4cm]{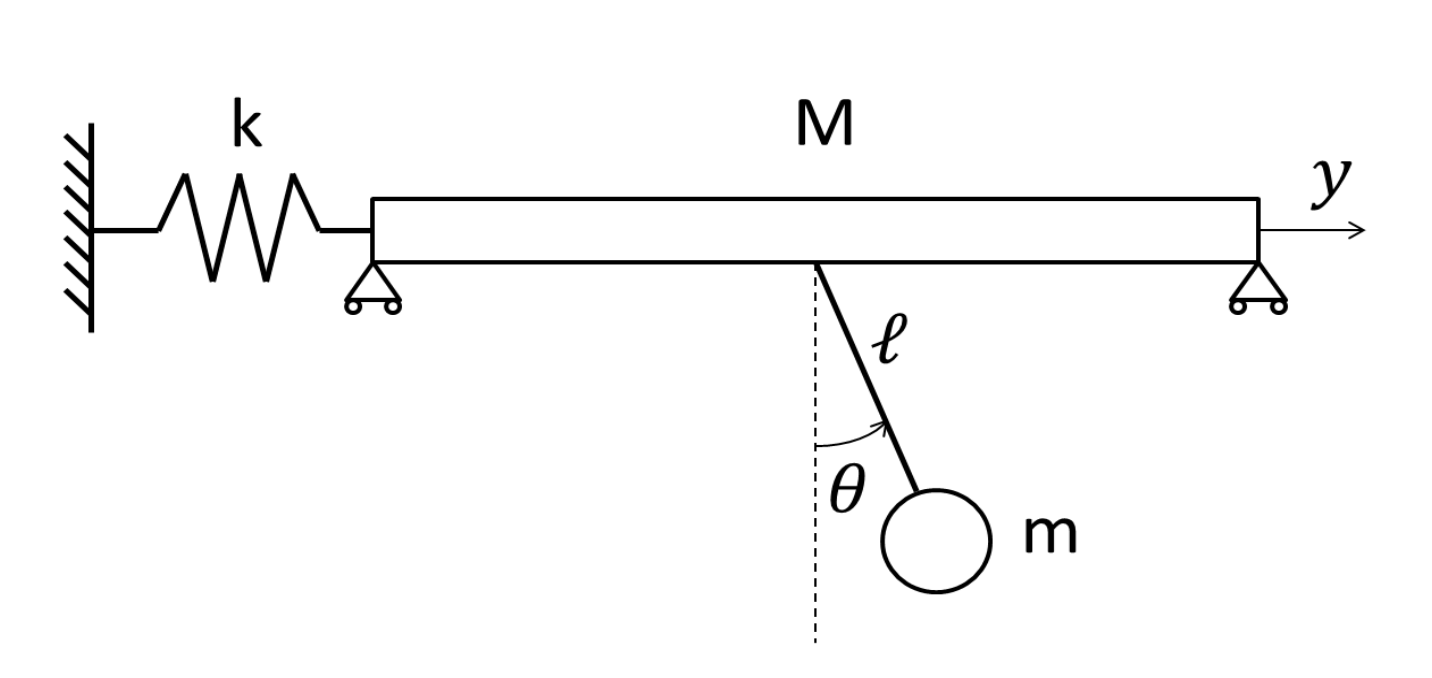}
    \caption{\label{fig::pendulumshema}The double pendulum system with its parameters.}
\end{figure}

\paragraph{The physical system}
The system can be seen as a dual-oscillator cluster.
The first oscillator is a spring-mass system whose axis is perpendicular to the gravitational axis.
The parameters of this system are the mass of the system $M_S$ and the spring stiffness $k$.
The initial position of the mass is denoted $y_0$, its initial velocity is $0$.
The second oscillator is a pendulum.
A schematic representation of the system is presented in \Cref{fig::pendulumshema}.
The parameters are the mass $m$ and the length of the pendulum $\ell$, which are fixed.
The initial value of the angle is $\theta_0$ and its derivative is $\dot{\theta_0}$.
By Newton's law of motion, the dynamics is governed by a system of two coupled ordinary differential equations (ODEs).
However, we do not have a closed form expression that gives the solution of the system.
This forces us to use computer codes.
The output signal is the position of the mass $m$ at time $t\in \left\{t_1,\ldots,t_{N_t}\right\}$ with $N_t=101$.
The input vector is $\bm{x}=\left\{ M_S, k, y_0,\theta_0,\dot{\theta_0}\right\}$.
The input variables are assumed to be independent and identically distributed with uniform distributions as described in \Cref{tab::inputs}.

\paragraph{The two different code levels}
We propose two codes.
The high-fidelity code numerically solves the coupled system of ODEs by an Euler's derivation of the position $y$ and the angle $\theta$ for each $t_u$.
This gives functions $\theta(t_u)$ and $y(t_u)$.
The low-fidelity code assumes that the angle of the $\theta$ pendulum is small so that the linearisation of $\sin({\theta})$ makes it possible to get a \textcolor{black}{simpler form of the expression of} the two coupled ODEs \textcolor{black}{and a faster resolution}.

\begin{table}[tbhp]
{\footnotesize
    \caption{\label{tab::inputs} Distributions of the input variables.}
    \centering
    \begin{center}
    \begin{tabular}{c | c | c | c | c}
        $M_S$ & $k$ & $\theta_0$ & $\dot{\theta_0}$ & $y_0$ \\ \hline
        $\mathcal{U}(10,12)$ & $\mathcal{U}(1;1.4)$ & $\mathcal{U}(\frac{\pi}{4};\frac{\pi}{3})$ & $\mathcal{U}(0;\frac{1}{10})$ & $\mathcal{U}(0;0.2)$ \\
    \end{tabular}
    \end{center}
}
\end{table}

\paragraph {Code Analysis}

\begin{figure}
    \centering
    \includegraphics[scale=0.3]{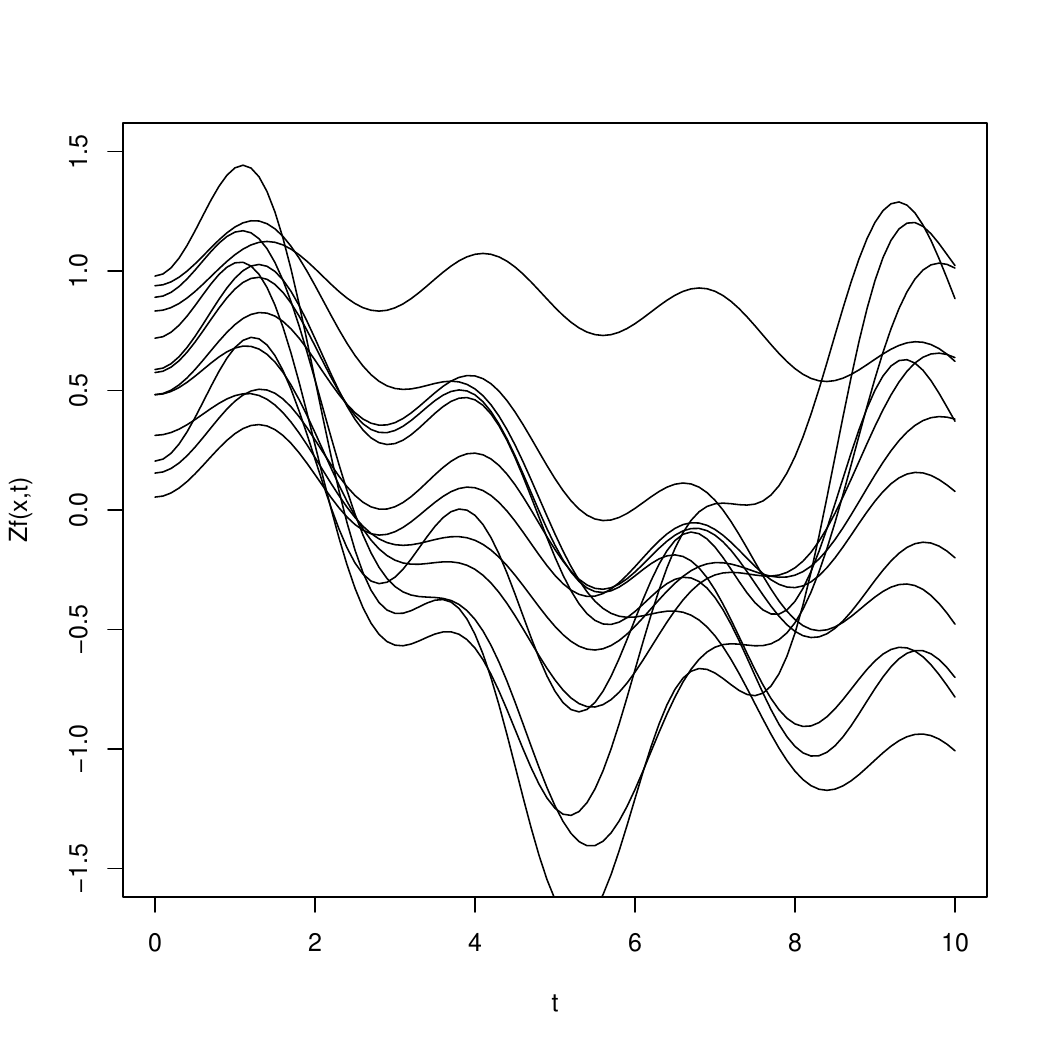}
    \includegraphics[scale=0.3]{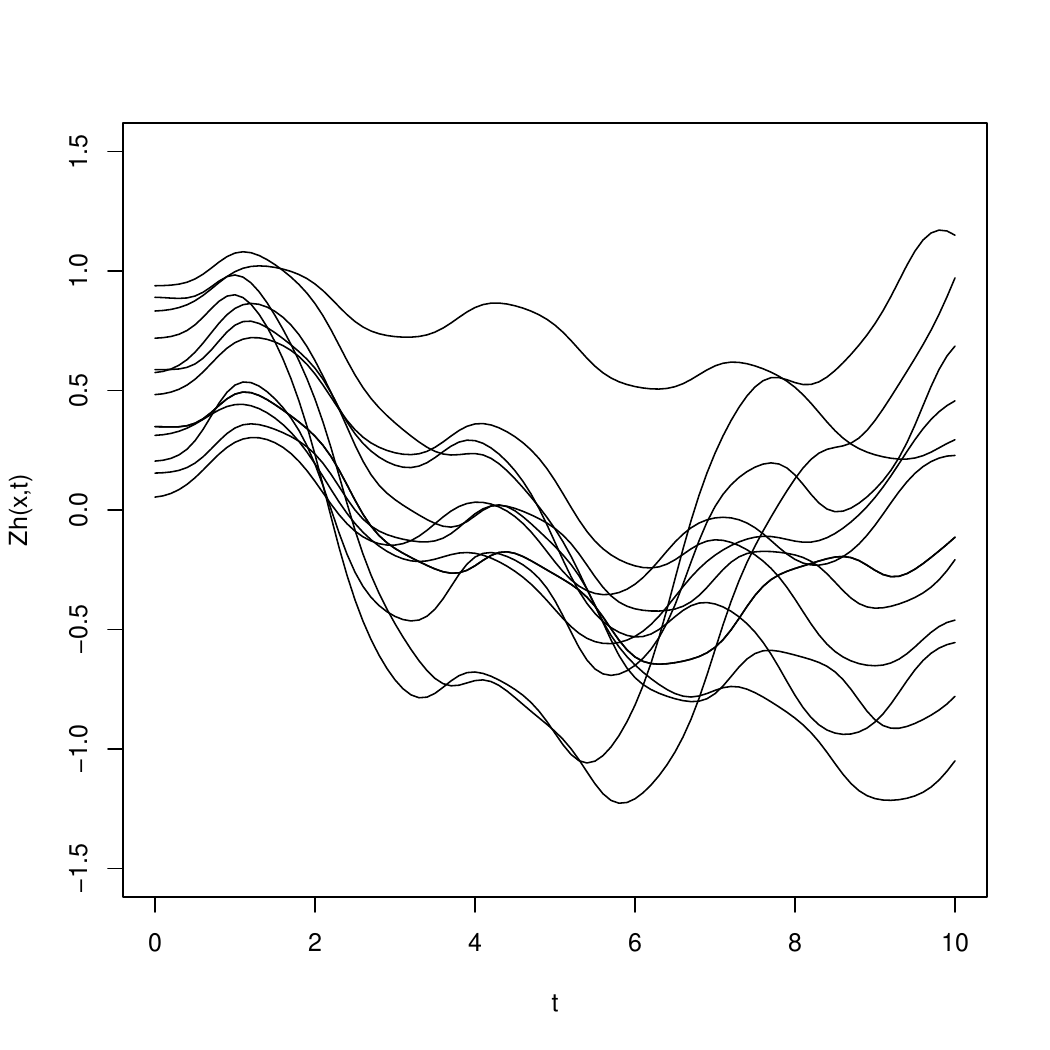}
    \includegraphics[scale=0.3]{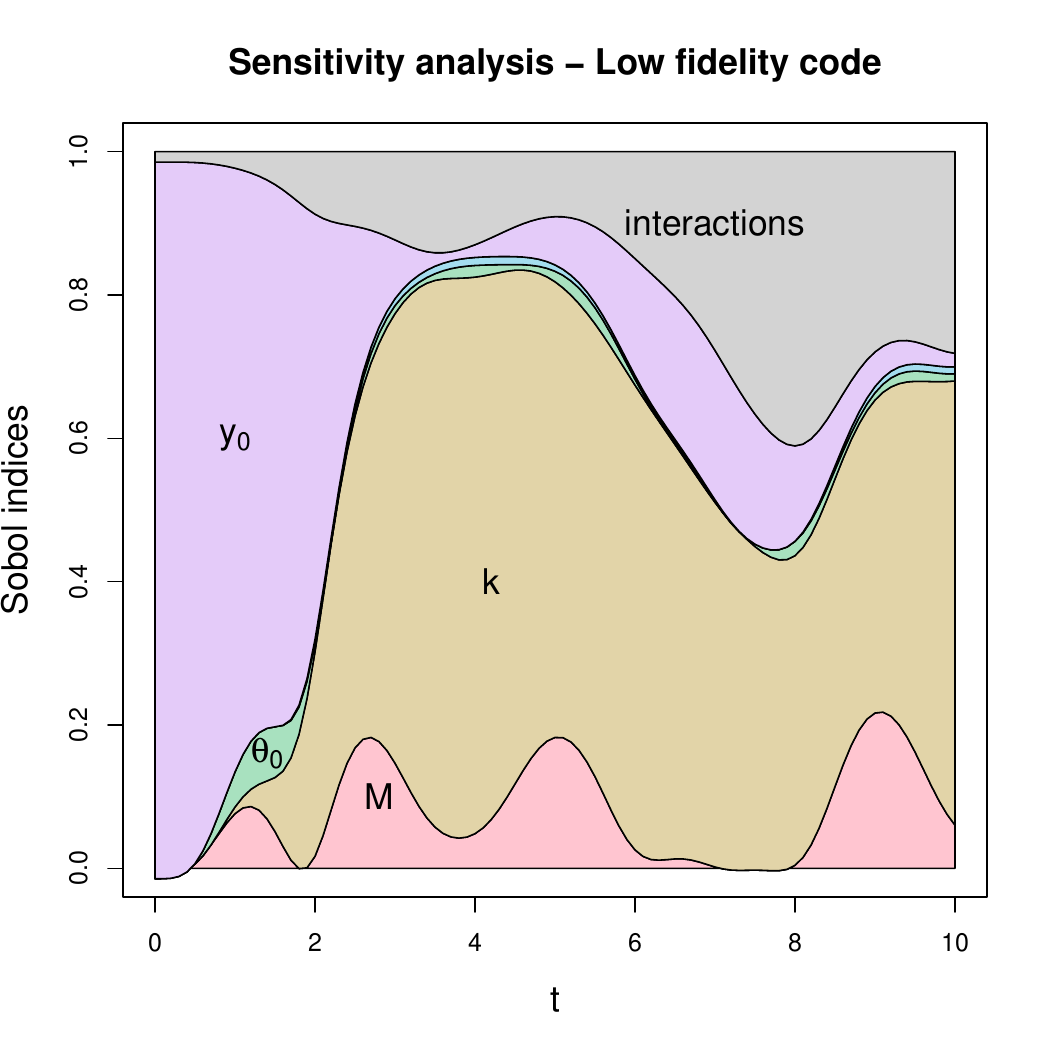}
    \includegraphics[scale=0.3]{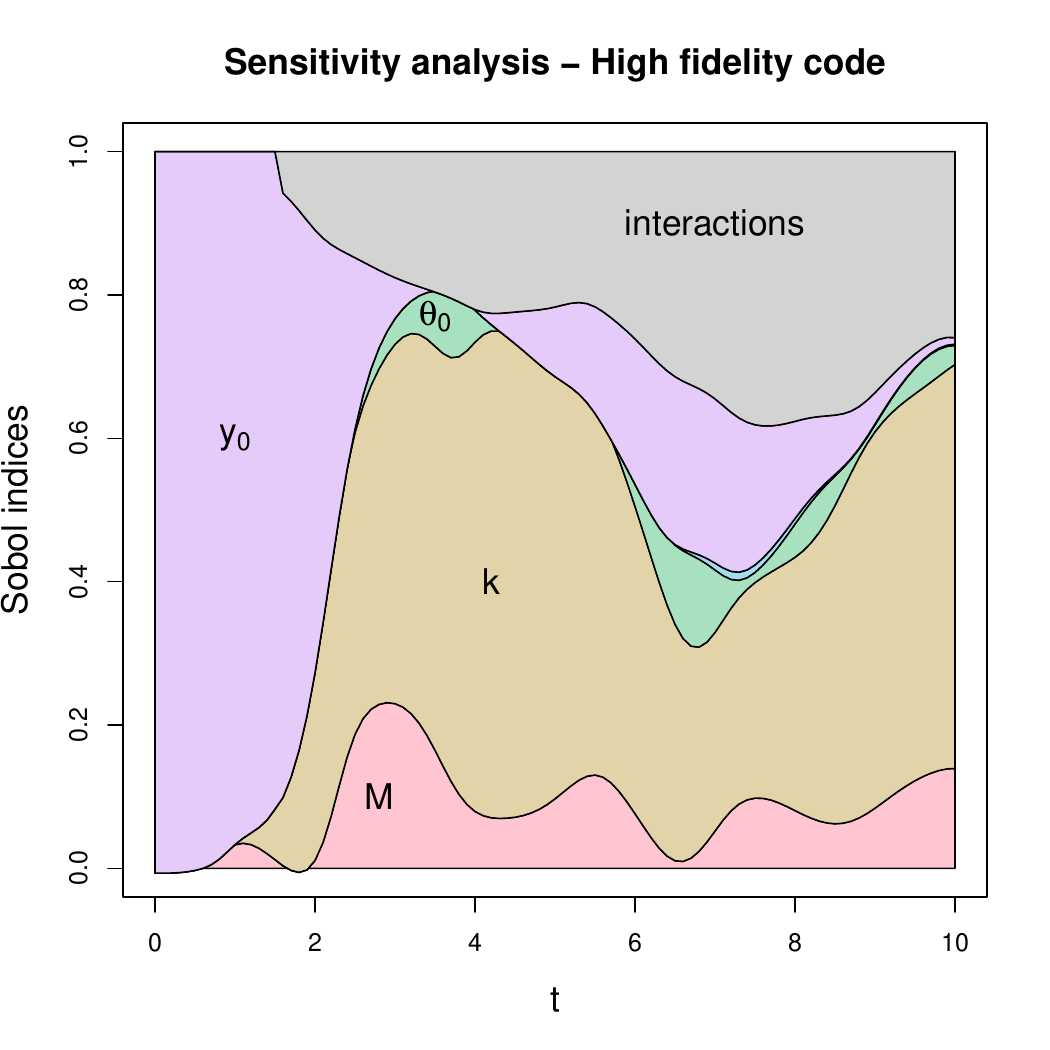}
    \includegraphics[scale=0.33]{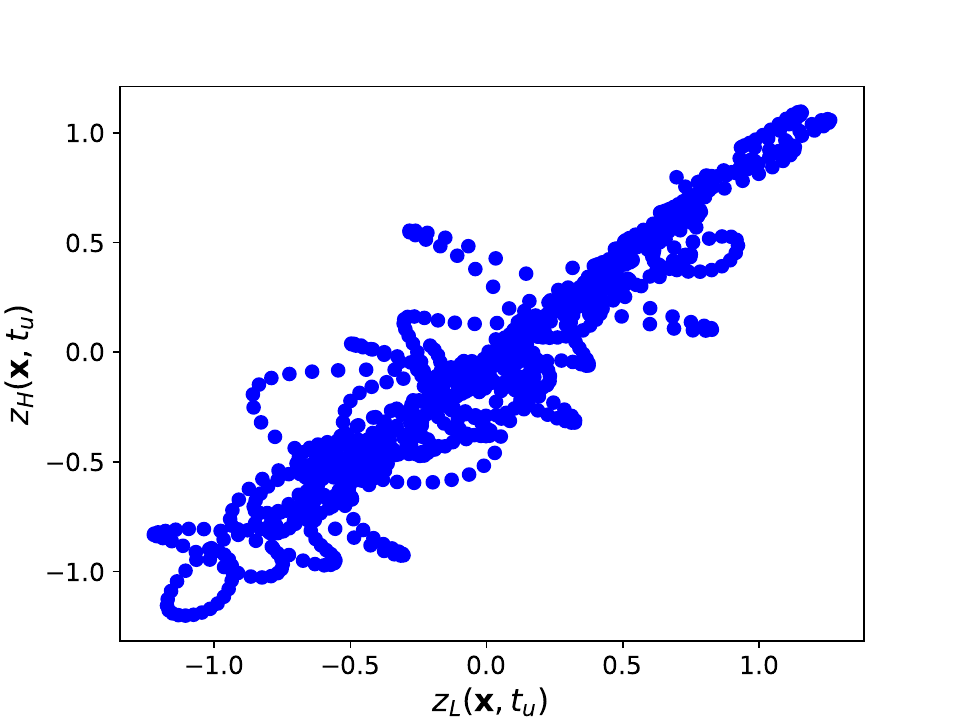}
    \includegraphics[scale=0.33]{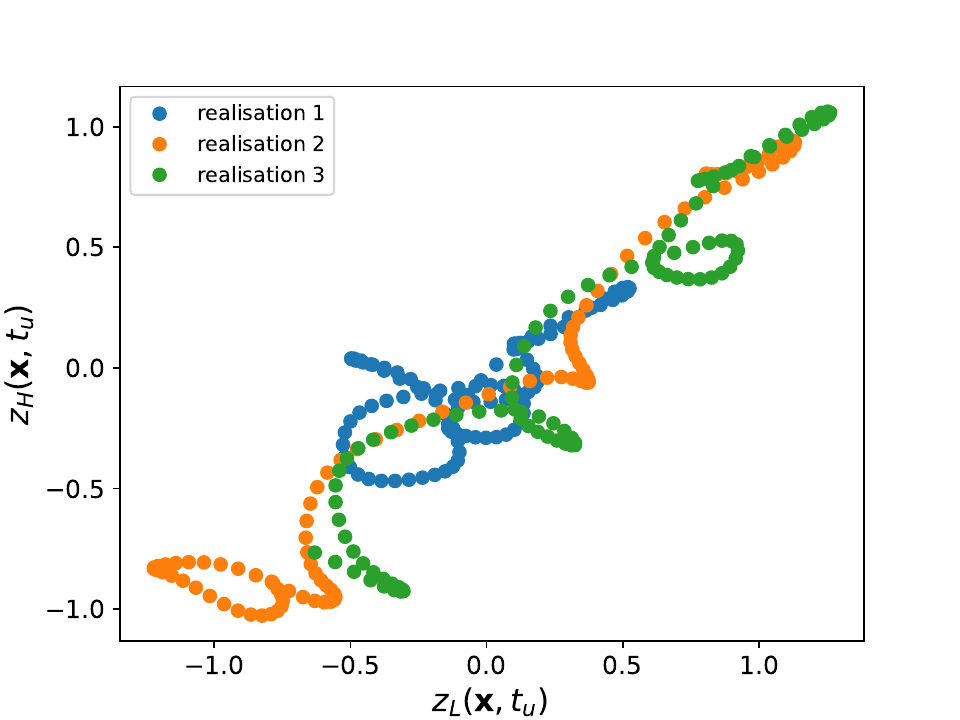}
    \caption{Comparison between low-fidelity (top left) and high-fidelity (top right) code outputs. 
        Sobol indices for high- and low-fidelity codes (\textcolor{black}{center} left: \textcolor{black}{low-}fidelity, \textcolor{black}{center} right: \textcolor{black}{high-}fidelity).
        For each time $t$ in the time grid, we report the first-order Sobol indices and "interactions" stands for the sum of the Sobol indices of order larger than $2$.
        \textcolor{black}{Finally, the bottom plots represent the interactions between codes (plots of $\left( z_L(x,t_u), z_H(x,t_u) \right)_{u=1}^{N_t}$ for different $x$).
        The bottom left graph is for 10 and the bottom right graph is for 3 values of $x$.}
        }
    \label{fig::SobolPend}
\end{figure}

A sensitivity analysis is carried out for information purposes, but it \textcolor{black}{is not used} in the forthcoming surrogate modeling.
\textcolor{black}{The sensitivity analysis makes it possible} to determine the effective dimension of our problem.
We compare outputs of the high- and low-fidelity codes and the associated Sobol indices on \Cref{fig::SobolPend}.
We estimate Sobol indices by the method described in \cite{saltelli2010variance} and implemented in the R library \cite{rSensitivity} by using a Monte Carlo sample of size $10^{5}$ for each code.
No surrogate model was used to estimate the indices in \Cref{fig::SobolPend}.
The main result is that the two codes depend on the same input variables.
The four most important input variables  are $y_0$, $k$, $M$ and $\theta_0$. 

\subsection{Comparison between methods}
The experimental designs used to compare the methods are presented in \cite{gratiet_multi-fidelity_2013}.
They are constructed from two independent maximum LHS designs with $N_\text{H}=10$ and $N_\text{L}=100$ points.
The low-fidelity design is then modified so that the designs are nested.
Only the points of the low-fidelity design closest to the points of the high fidelity design are moved.
To generate these designs the R packages \cite{rDiceDesEval,rMuFiCokriging} \textcolor{black}{are} used.
\textcolor{black}{A random uniform nested design can also be used, but we choose a more effective design for GP regression.}
The test design is composed of $4000$ \textcolor{black}{points} randomly chosen in the hypercube determined by the supports of the uniform distributions described in \Cref{tab::inputs}.

In this section, we want to demonstrate the interest of the method presented in \Cref{sec::AR1multiTens}.
For this we will compare several methods:
\begin{itemize}
    \item the multi-fidelity method presented in \Cref{sec::AR1multiproj} with a Dirac distribution of $\bm{\Gamma}$, called the \textcolor{black}{SVD} method.
    \item the multi-fidelity method that uses \textcolor{black}{GP} regression of the orthogonal part \textcolor{black}{with covariance-tensorization} and the distribution of $\bm{\Gamma}$ is Dirac at $\bm{\gamma}$ the matrix of the SVD of the observed low-fidelity code outputs.
    It prediction is computed as is \Cref{sec::AR1multiTens} and called Dirac method. 
    \item the multi-fidelity method presented in \Cref{sec::AR1multiTens} with the \textcolor{black}{CVB} distribution, called \textcolor{black}{CVB} method.
    \textcolor{black}{\item the neural network (NN) method presented in \cite{meng2020composite}.
    We extend this method for time-series outputs by considering $N_t$-dimensional outputs for the low- and high-fidelity neural networks and by removing the physical inspired NN part.
    We used the parameters proposed in the article, i.e. $2$ hidden layers for each network with $20$ neurons per layer.
    We also tested the NN method up to $100$ neurons per layer but the best results were obtained with $20$ neurons approximately.}
\end{itemize}
The method we would like to highlight is the \textcolor{black}{CVB} method.

\paragraph{\textcolor{black}{CVB} basis}
The law of $\bm{\Gamma}$ needs to be determined in order to compute or estimate the moments (\ref{eq::MeanZpara}),(\ref{eq::VarZparaFull}), (\ref{eq::meanArbTens}), and (\ref{eq::VarFullortho}).
The distribution of $\bm{\Gamma}$ is the \textcolor{black}{CVB} distribution described in \Cref{sssec::ELprojbase}.
As shown by \cref{eq::emplowGamma}
it depends on the size $k$ of the random subset $I$.
Here we choose $k=4$.
Because it is too expensive to compute the sum over all $\binom{N_\text{L}}{k}$ different subsets $\{j_1,\hdots,j_k\}$,
we estimate the expectation (\ref{eq::emplowGamma}) by an empirical average over $n=64$ realizations $I_j$ of the random subset $I$:
\begin{equation}
    \label{eq::emplowGammaemp}
    \E\left[f(\bm{\Gamma})\right] \simeq \frac{1}{n}
    \sum_{j=1}^n f(\bm{\tilde{U}}_{I_j}) ,
\end{equation}
where $\bm{\tilde{U}}_{I_j}$ is the matrix of the left  singular vectors of the SVD of\\ $\left(z_\text{L}(\bm{x}^{(i)},t_u)\right)_{\substack{u\in\{1,\ldots,N_t\}\\ i\in \{1,\ldots,N_\text{L}\}\backslash I_j}}$.
We have checked that the stability with respect to $k$ is conserved if $1<k<N_\text{L}-N_\text{H}$ and the stability with regard to $n$ is valid if $n>\max{(k, 50)}$.
We \textcolor{black}{have} tested the construction of the \textcolor{black}{CVB} basis for all $k$ values in this range and found that changes in $k$ do not influence the basis \textcolor{black}{significantly}.

The computational cost of calculating the basis is very important in particular because it is impossible for us to calculate it for all subsets.
A method to compute the basis with only a cost of $O(N_t^2)$ is given in \cite{mertens1995efficient} whereas we compute it with $O(N_\text{L}^2N_t)$ by our method.
The gain is however very small especially if $N_\text{L}\ll N_t$ which is our case.
We have therefore not implemented this method in the results presented \textcolor{black}{in this paper}.

\paragraph{Prediction of the orthogonal part}
A simple model for the a priori mean function is \textcolor{black}{chosen} $M=1$ and $f(\bm{x})=1$.
Consequently, $F^TR_x^{-1}F=\sum_{i,j}\left\{R_x^{-1}\right\}_{i,j}$.

\paragraph{Multi-fidelity regression of the coefficients}

Our implementation of the multi-fidelity regression is based on \cite{rMuFiCokriging}.
We use an informative prior for the regression of the coefficients.
For more information refer to \cite[Section 3.4.4]{gratiet_multi-fidelity_2013}.
In this example the size of the p\textcolor{black}{r}iors are $q=p_\text{L}=p_\text{H}=1$.
Considering the relation between the two codes we cho\textcolor{black}{o}se $b^\rho = 1$.
The trend i\textcolor{black}{s} supposed to be null consequently, $b^\beta_\text{H} = b_\text{L} = 0$.
The variances are $\sigma_\text{L} = 0.5$ and $\sigma_\text{H}=0.5$ with $V_\text{H}^\beta=2$ and $V_\text{L}=2$.
The parameters for the inverse Gamma distribution $m_\text{L}= m_\text{H} = 0.2$ and $\varsigma_\text{L} = \varsigma_\text{H} = 1.5$.
We have checked the robustness of the results with respect to the hyper-parameters of the prior distributions.
Alternatively, the article \cite{ma2020objective} presents non-informative priors for the autoregressive co-kriging.

\paragraph{Prediction}

\begin{figure}
    \centering
    \includegraphics [scale=0.5]{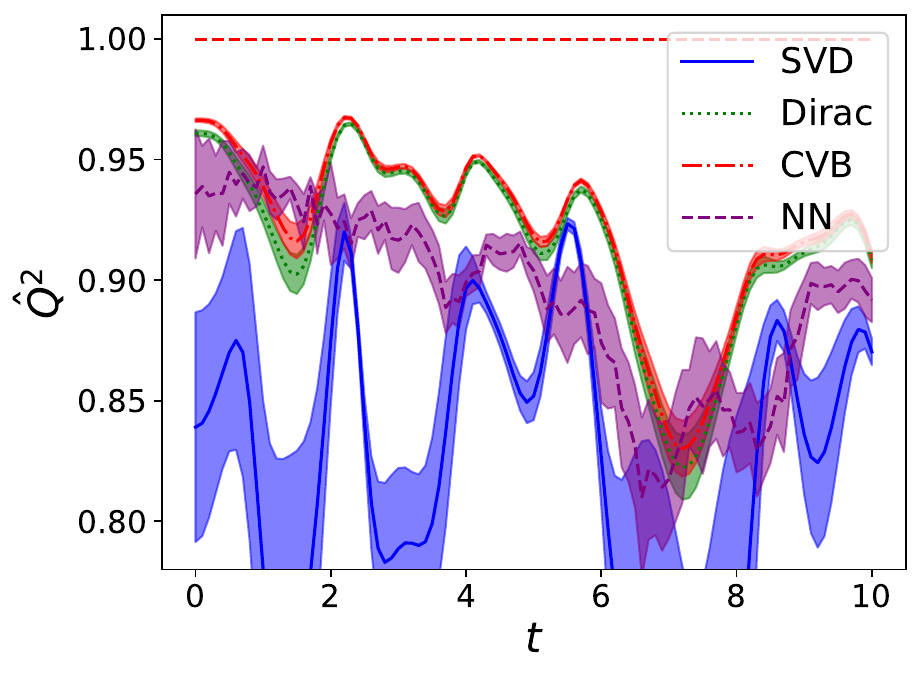}
    \caption{\label{fig::ResultQ2tot}Comparison between the methods in terms of time-dependent $\hat{Q}^2$.
        Averages over $40$ random experimental designs are computed.
    The  colored fields represent the confidence intervals determined by $\pm1.96$ empirical standard deviation. Here
    $N_\text{H}=10$, $N_\text{L}=100$ and $N_t=101$.}
\end{figure}

In order to estimate the errors of the surrogate models, we calculate their $\hat{Q}^2$'s and report them in \Cref{fig::ResultQ2tot}.
To compute $\hat{Q}^2$ for a model, we calculate the difference between \textcolor{black}{the} validation set of size $4000$ and the predictions of the model.
We have averaged the estimates of the $\hat{Q}^2$ over $40$ different experimental designs.

\begin{figure}
    
    \centering
    \includegraphics[scale=0.33]{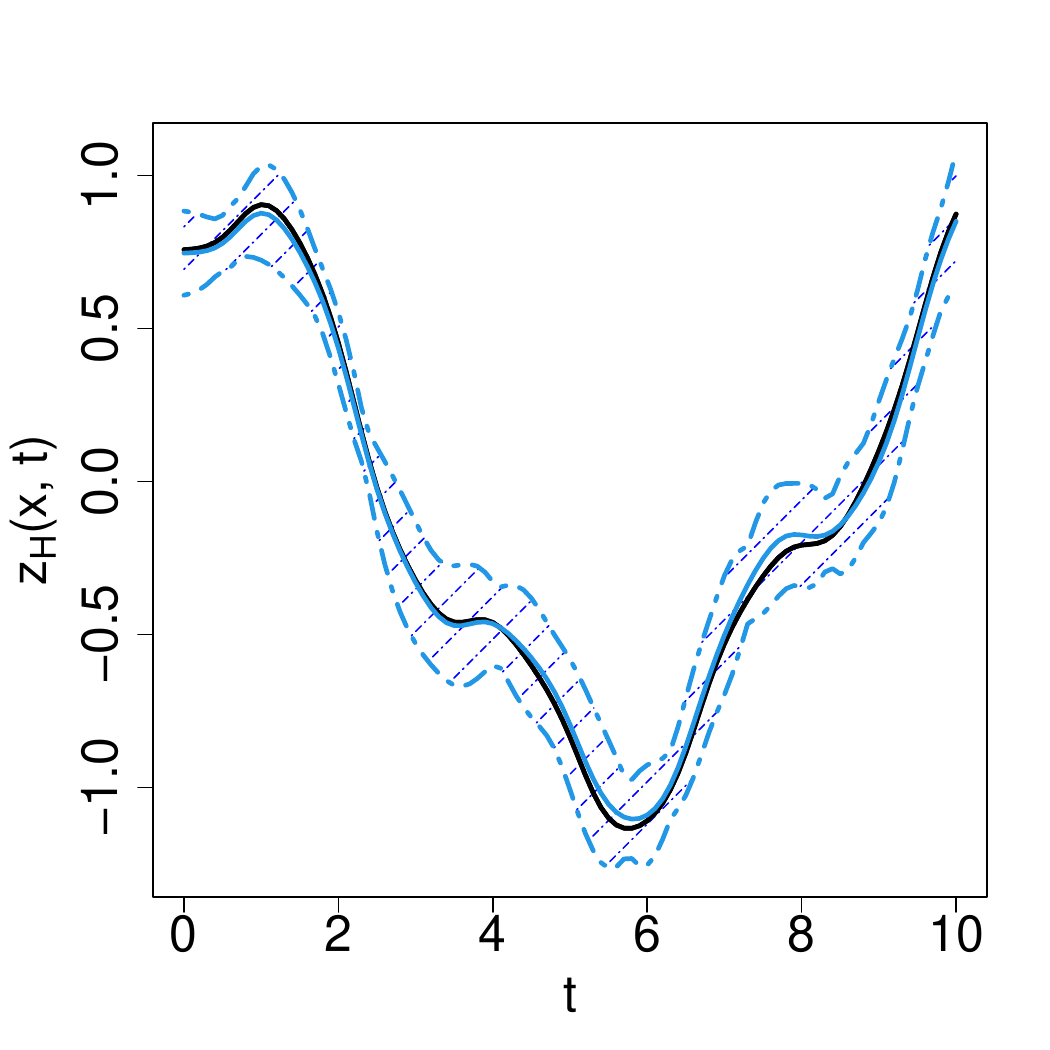}
    \includegraphics[scale=0.33]{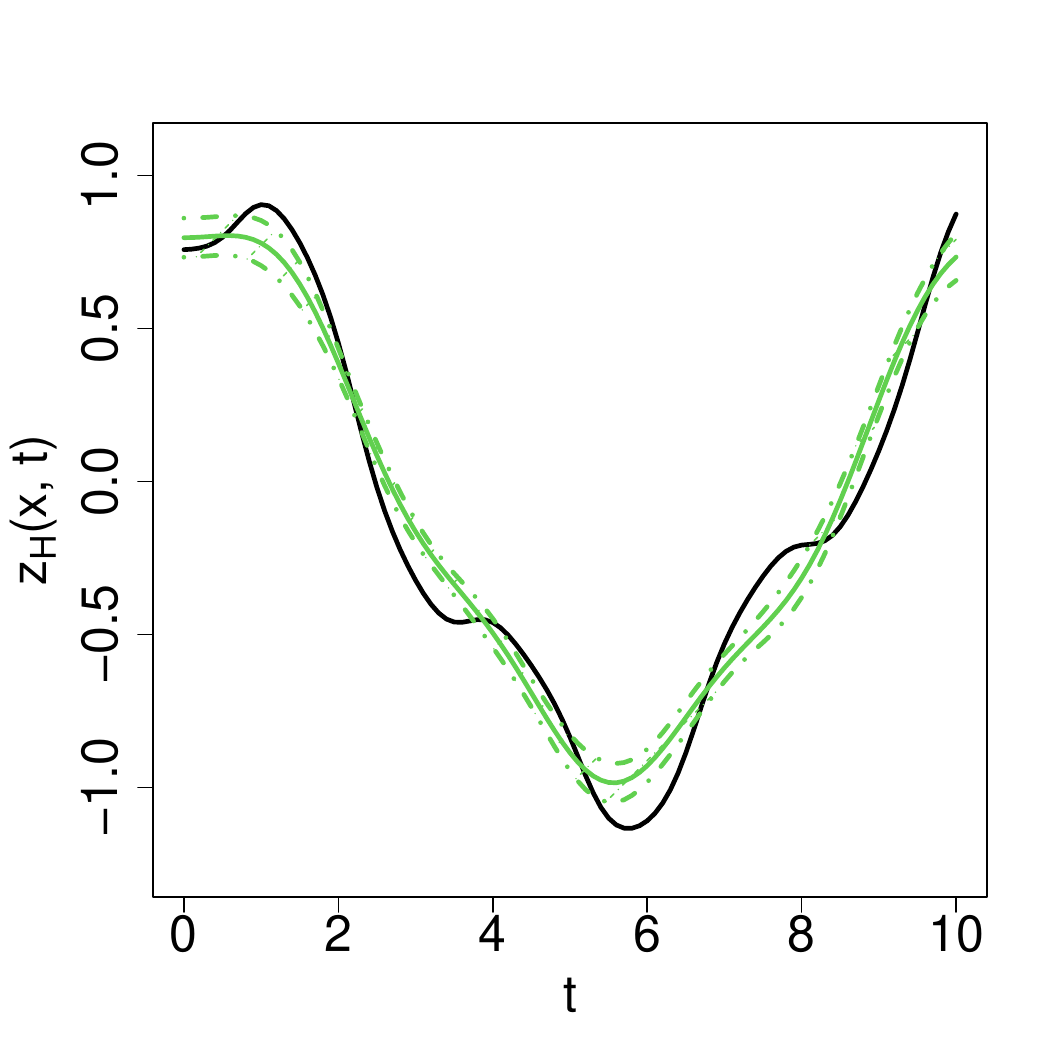}
    \caption{\label{fig::exampleresultcurve} Comparison between the predictions of the \textcolor{black}{CVB} method (left) and the \textcolor{black}{SVD} method (right).
    The black solid line is the exact high-fidelity time series, the colored solid line is the prediction mean and the dashed lines are the confidence intervals.
    In this example the value of $N$ obtained by cross validation is $8$.}
\end{figure}

The \textcolor{black}{SVD} method gives a very interesting result because the $\hat{Q}^2$ is almost always higher than $0.8$.
However, in \Cref{fig::exampleresultcurve} we can see that it does not capture the form of the times series.
The $\hat{Q}^2$ of the Dirac and \textcolor{black}{CVB} methods are larger 
 than the ones of the other methods.
The error is also less variable as a function of $t$.
And the variance is much lower for both methods.
However, even if there is a difference between the Dirac and \textcolor{black}{CVB} methods, it is not possible to say that the \textcolor{black}{CVB} method is better in this application.
The difference between the Dirac and \textcolor{black}{CVB} methods is small, in our example.

The variance of the prediction is very important for the quantification of \textcolor{black}{prediction} uncertainty.
All formulas are given in  the previous sections and we illustrate the results in \Cref{fig::exampleresultcurve}.
We can see that the variance of the projection method is not accurate and overestimates the quality of the prediction.
This method is not acceptable for prediction.
The Dirac method and the \textcolor{black}{CVB} method have almost the same variance.
If we compare to the variance of the \textcolor{black}{SVD} method, it means that most of the uncertainty relates to the orthogonal part.
This leads to the conclusion that this part is important in the regression.

\begin{figure}
    \label{fig::variance}
    \centering
    \includegraphics[scale=0.4]{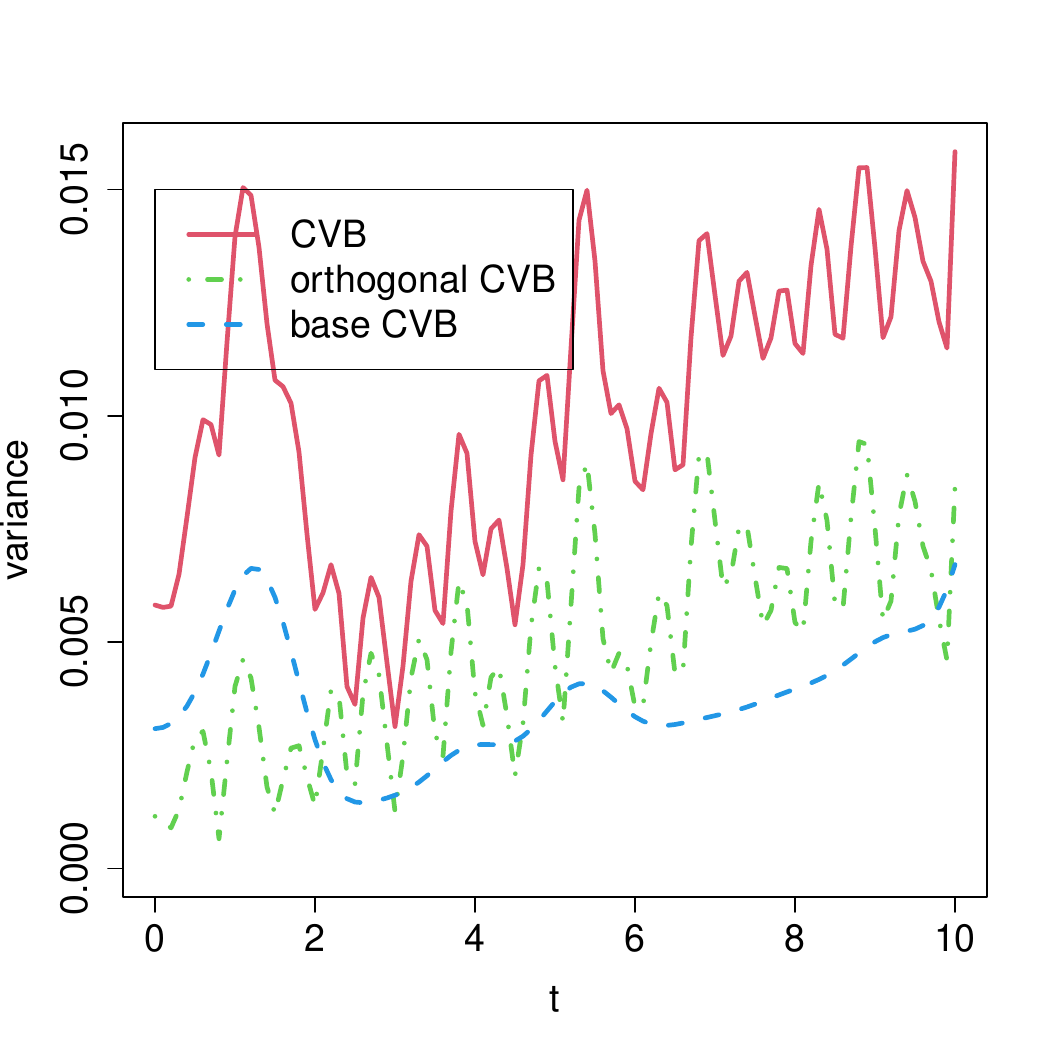}
    \caption{Estimation of the different time-dependent prediction variance terms for the empirical method.}
\end{figure}

In order to understand the interest of the method with covariance tensorization for the orthogonal part we study in more detail the orthogonal part.
First we study \textcolor{black}{the role of} the value of $N$.
Here $N_\text{L}=100$, \textcolor{black}{so the possible values of $N$ are between $0$ and $100$}.
We find that the optimal value of $N$ for $40$ learning sets is between $8$ and $10$.
Even when the value of $N_\text{H}$ is increased, $N$ remains constant in the $8$ to $10$ range.
This means that the low-fidelity code can give reliable information on the high-fidelity code output  projection into a $8$-dimensional space. The high-fidelity code output is, however, higher dimensional and it is important to predict the orthogonal part with a dedicated method, namely the proposed covariance tensorization method.

\textcolor{black}{
We have carried out extensive numerical simulations with values of $N_H$ in the range $\left[ 5:20 \right]$ and values of $N_L$ in the range $\left[50:1000\right]$.
If only very small data sets are available ($5 \leq N_H \leq 7$) the prediction is not satisfactory whatever the method.
Moreover, for values of $N_L$ greater than 200, there is no significant change except for the NN method which improves to the level of the CVB method for the largest data sets.
There is also a decrease in prediction uncertainty with the increase of the $N_H$ number as can be expected.
At the same time as the prediction uncertainty increases when $N_H$ is decreased, there is also a decrease in prediction performance, but independently of the method.
The code we used is available in \cite{Kerleguer_MultiFi_Time-Series_2022}.
}



%

\section{Discussion}
The objective of this work is to propose a method that generates a surrogate model in the context of multi-fidelity and time series outputs \textcolor{black}{and that} quantific\textcolor{black}{es the prediction} uncertainty.
The method we propose is based on three main ingredients: dimension reduction, co-kriging, and covariance tensorization.
The model we present is based on multi-fidelity (co-kriging) regression.
By reducing the output dimension, multi-fidelity regression becomes possible.
To take into account all the information \textcolor{black}{contained in} the data sets, the part that cannot be treated with the previous method is predicted by Gaussian process regression with covariance tensorization.

First, we have presented different ways to build the basis that allows  to represent the high\textcolor{black}{-} and low\textcolor{black}{-}fidelity code outputs.
Second, we have presented a model that allows to estimate the high-fidelity code outputs from data collected from the high- and low-fidelity codes.
The combination of a multi-fidelity part and a single-fidelity part with tensorized covariance is the central point of the proposed method.
The performance of our model has been tested on a mechanical example.
We have been able to use multi-fidelity in a very convincing way to build a robust surrogate model better than any other method presented so far.

There are several ways to extend the method presented in this article.
Sequential experimental designs in a multi-fidelity context have already been dealt with by \cite{le2015cokriging}.
However, they deserve to be extended to the case of time-series outputs.
We can consider regression problems for more than two levels of code.
It is conceivable in this case to build several levels of bases which from code to code would improve the basis and thus reduce the orthogonal part.
In addition, high-dimensional outputs are not different from time-series outputs as considered in this paper.
It is therefore conceivable to adapt this method to more general functional outputs.

\appendix

\section{Multi-fidelity priors for AR(1) model}
\label{asec::MultiFi}
\textcolor{black}{In the following we define the priors needed to use the AR(1) model defined in \cref{subsec::multiGPreg}.}

The goal of a Bayesian prediction is to integrate the uncertainty of the parameter estimation into the predictive distribution as in \cite{le2013bayesian}.
Here the parameters are $\sigma$, $\beta$ and $\beta_\rho$.
As explained in \cite{le2014recursive} the result is not Gaussian but we can obtain expressions of the posterior mean $\E\left[A_\text{H}(\bm{x})|\mathcal{A}=\alpha\right]$ and variance $\V\left[A_\text{H}(\bm{x})|\mathcal{A}=\alpha\right]$.
It is possible to consider informative or non informative priors for the parameters \cite{le2014recursive,ma2020objective}.
Here we consider informative conjugate priors:
\begin{align}
        \left[\sigma^2_\text{L}\right] \sim & \mathcal{IG}\left(m_\text{L},\varsigma_\text{L}\right), \\
        \left[\beta_\text{L}|\sigma_\text{L}^2\right]
        \sim & \mathcal{N}_{p_\text{L}} \left(b_\text{L},\sigma^2_\text{L}V_\text{L}\right),\\
        \left[\sigma^2_\text{H}|\mathcal{A}^{\text{L}}=\alpha^{\text{L}},\beta_\text{L},\sigma_\text{L}\right] \sim &\mathcal{IG}\left(m_\text{H},\varsigma_\text{H}\right),\\
        \left[\beta_\rho, \beta_\text{H}|\mathcal{A}^{\text{L}}=\alpha^{\text{L}},\sigma_\text{H}^2, \beta_\text{L}, \sigma_\text{L}\right]
        \sim & \mathcal{N}_{q+p_\text{H}}\left(b_\text{H}=\begin{pmatrix}b^\rho\\ b_\text{H}^\beta\end{pmatrix}, \sigma_\text{H}^2 V_\text{H} = \sigma_\text{H}^2\begin{pmatrix}V^\rho & 0\\ 0 & V_\text{H}^\beta\end{pmatrix}\right),
    \end{align}
with \begin{itemize}[noitemsep,topsep=0pt,parsep=0pt,partopsep=0pt]
    \item[-] $b_\text{L}$ a vector of size $p_\text{L}$,
    \item[-] $b^\rho$ a vector of size $q$,
    \item[-] $b^\beta_\text{H}$ a vector of size $p_\text{H}$,
    \item[-] $V_\text{H}^\beta$ a $p_\text{H}\times p_\text{H}$ matrix,
    \item[-] $V^\rho$ a $q\times q$ matrix,
    \item[-] $V_\text{L}$ a $p_\text{L}\times p_\text{L}$ matrix,
    \item[-] $m_\text{F}$ and $\varsigma_\text{F}$ are positive \textcolor{black}{real numbers} and $\mathcal{IG}$ stands for the inverse Gamma distribution.
\end{itemize} 
By using these informative conjugate priors we obtain the following a posteriori distribution\textcolor{black}{s} as in \cite{le2014recursive} :

\begin{align}
        \left[\sigma_\text{L}^2 | \mathcal{A}^\text{L}=\alpha^\text{L} \right] \sim  & \mathcal{IG}\left(d_\text{L}, Q_\text{L}\right) , \\
        \left[\beta_\text{L} | \mathcal{A}^\text{L}=\alpha^\text{L}, \sigma_\text{L}^2 \right] \sim & \mathcal{N}_{p_\text{L}} \left(\Sigma_\text{L} \nu_\text{L}, \Sigma_\text{L}\right) , \\
        \left[\sigma_\text{H}^2 | \mathcal{A}=\alpha \right] \sim & \mathcal{IG}\left(d_\text{H}, Q_\text{H}\right) , \\
        \left[\beta_\text{H}, \beta_{\rho} | \mathcal{A}=\alpha, \sigma_\text{H}^2\right] \sim & \mathcal{N}_{p_\text{H}+q}\left(\Sigma_\text{H}\nu_\text{H},\Sigma_\text{H}\right) ,
     \label{eq::priorMultiLoic}
\end{align}
with:
\begin{itemize}[noitemsep,topsep=0pt,parsep=0pt,partopsep=0pt]
    \item [-] $d_\text{F}=\frac{n_\text{F}}{2}+m_\text{F},$
    \item [-] $\tilde{Q}_\text{F}=\left(\alpha^{\text{F}}-H_\text{F}\hat{\lambda}_\text{F}\right)^T C_\text{F}^{-1}\left(\alpha^{\text{F}}-H_\text{F}\hat{\lambda}_\text{F}\right),$
    \item [-] $Q_\text{F} = \tilde{Q}_\text{F} + \varsigma_\text{F} + \left(b_\text{F}-\hat{\lambda}_\text{F}\right)^T\left(V_\text{F}+\left(H_\text{F}^TC_\text{F}^{-1}H_\text{F}\right)\right)^{-1}\left(b_\text{F}-\hat{\lambda}_\text{F}\right)$,
    \item [-] $\Sigma_\text{F} = \left[H_\text{F}^T\frac{C_\text{F}^{-1}}{\sigma^2_\text{F}}H_\text{F} + \frac{V_\text{F}^{-1}}{\sigma^2_\text{F}}\right]^{-1}$,
    \item [-] $\nu_\text{F} = \left[ H_\text{F}^T\frac{C_\text{F}^{-1}}{\sigma_\text{F}^2}\alpha^{\text{F}}+ \frac{V_\text{F}^{-1}}{\sigma^2_\text{F}}b_\text{F}\right]$,
    \item [-] $H_\text{F}$ is defined by $H_\text{L}=F_\text{L}$ and $H_\text{H} = \left[G^\text{L}\odot\left(\alpha^{\text{H}}\bm{1}^T_{q_\text{L}} \right)~F_\text{H}\right] $,
    \item [-] $G^\text{L}$ is the $N_\text{H}\times q$ matrix containing the values of $g_\text{L}^T(\bm{x})$ for $\bm{x} \in D_\text{H}$,
    \item [-] $\bm{1}_{q_\text{L}}$ is a $q$-dimensional vector containing $1$,
    \item [-] $\hat{\lambda}_\text{F} = \left(H_\text{F}^T C_\text{F}^{-1}H_\text{F} \right)^{-1} H_\text{F}^T C_\text{F}^{-1} \alpha^{\text{F}}$.
\end{itemize}

Consequently, the posterior distribution of $A_\text{H}(\bm{x})$ has the following mean and variance:
\begin{align}
    \label{eq::FullUKMultiFi}
    \E\left[A_\text{H}(\bm{x})|\mathcal{A}=\alpha\right] = &
    h_\text{H}^T(\bm{x})\Sigma_\text{H}\nu_\text{H}
    + r_\text{H}^T(\bm{x})C_\text{H}^{-1}\left(\alpha^{\text{H}}-H_\text{H}\Sigma_\text{H}\nu_\text{H}\right),\\
\nonumber
        \V\left[A_\text{H}(\bm{x})|\mathcal{A}=\alpha\right] =&
        \left(\hat{\rho}^2_\text{L}(\bm{x}) + \varepsilon_\rho (\bm{x})\right) \sigma_{\tilde{A}_\text{L}}^2(\bm{x})
        + \frac{Q_\text{H}}{2(d_\text{H}-1)}(1-r_\text{H}^T(\bm{x})C_\text{H}^{-1}r_\text{H}(\bm{x}))\\
        &+ \left(h_\text{H}^T-r_\text{H}^T(\bm{x})C_\text{H}^{-1}H_\text{H}\right)\Sigma_\text{H}\left(h_\text{H}^T -r_\text{H}^T(\bm{x})C_\text{H}^{-1}H_\text{H}\right)^T  ,
    \label{eq::FullUKMultiVar}
\end{align}
with $\hat{\rho}_\text{L}(\bm{x})= g_\text{L}^T(\bm{x})\hat{\beta_\rho}$,
$\hat{\beta}_\rho=\left[\Sigma_\text{H}\nu_\text{H}\right]_{i=p_\text{H}+1,\ldots,p_\text{H}+q}$ and 
$\varepsilon_\rho (\bm{x})= g_\text{L}^T(\bm{x})\tilde{\Sigma}_\text{H}g_\text{L}(\bm{x})$ with $\tilde{\Sigma}_\text{H} = \left[\Sigma_\text{H}\right]_{i,j =p_\text{H}+1,\ldots,p_\text{H}+q}$.

The posterior mean $\E\left[A_\text{H}(\bm{x})|\mathcal{A}=\alpha\right]$ is the predictive model of the high fidelity response and the posterior variance $\V\left[A_\text{H}(\bm{x})|\mathcal{A}=\alpha\right]$ represents the predictive variance of the model.

\section{LOO formula and discussion}
\label{asec::LOO}
\paragraph{LOO without loop}
In order to quickly minimize the LOO error with respect to the vector of correlation lengths $\bm{\ell_x}$ there exist formulas to evaluate $\varepsilon^2(\bm{\ell_x})$ with matrix products \cite{bachoc2013cross},\cite{dubrule1983cross}.
The LOO optimization problem is equivalent to minimize a function $f_{\textbf{CV}}(\bm{\ell_x})$ given by:
\begin{equation}\label{LOOcostfunct}
f_{\textbf{CV}}(\bm{\ell_x}) = z^T R_{\bm{\ell_x}}^{-1} \diag \left(R_{\bm{\ell_x}}^{-1}\right)^{-2} R_{\bm{\ell_x}}^{-1} z,
\end{equation} 
where $z$ is a vector that collect all the data.

Considering \Cref{matcorsep} and the mixed-product property, the inverse of a Kronecker product  and the formula $\diag(A \otimes B)=\diag A \otimes \diag B$ the cost function can be expressed as:
\begin{equation} \label{LOOcostfunctred}
f_{\textbf{CV}}(\bm{\ell_x}) = z^T \left(R_{t}^{-1} \diag \left(R_{t}^{-1}\right)^{-2} R_{t}^{-1}\right) \otimes \left(R_{x}^{-1} \diag \left(R_{x}^{-1}\right)^{-2} R_{x}^{-1}\right) z.
\end{equation}

 However, the term in $R_{t}$ is impossible to calculate because $R_{t}$ is not invertible.
$R_{t}^{-1} \diag \left(R_{t}^{-1}\right)^{-2} R_{t}^{-1}$ can be approximated by $I_{N_t}$ in order to have a tractable problem.
This assertion is equivalent to the hypothesis:
\begin{equation}\label{conditionmatcor}
R_t^{2}=\diag(R_t^{-1})^{-2} 
.
\end{equation}
This assumption can be seen as the fact that the error is estimated by taking into account only the spatial distribution of the covariance.
Indeed, to calculate the error only the matrix $R_x$ is used, even if the value of $R_t$ is calculated by maximum likelihood later in the \textcolor{black}{GP} regression method.

Thus, the minimization described in \Cref{eq::LOOerrnorm} makes it possible to calculate the correlation lengths by minimizing:

\begin{equation} \label{LOOcostsimp}
f_{\textbf{CV}}(\bm{\ell_x}) \simeq z^T I_{N_t} \otimes \left(R_{x}^{-1} \diag \left(R_{x}^{-1}\right)^{-2} R_{x}^{-1}\right) z ,
\end{equation}
\textcolor{black}{where $I_{N_t}$ is the $N_t\times N_t$ identity matrix.}
The main interest of this method is to give an approximate value of the error and to make the optimization much faster.

\paragraph{Optimization with hypothesis (\ref{conditionmatcor})} Efficient minimization algorithms require  to have the derivative of the function so that it does not have to be calculated by finite differences.
Thanks to the simplification (\ref{LOOcostsimp}) it is possible to calculate the derivative of the LOO error \cite{bachoc2013cross}:

\begin{equation}\label{dereqsimp}
\begin{aligned}
\frac{\partial}{\partial \ell_{x_k}}f_{\textbf{CV}}(\bm{\ell_x}) = 2 z^T I_{N_t} \otimes R_{x}^{-1} \diag \left(R_{x}^{-1}\right)^{-2} \left( R_{x}^{-1} \frac{\partial R_{x}}{\partial \ell_{x_k}} R_{x}^{-1} \right) \diag \left(R_{x}^{-1}\right)^{-1} R_{x}^{-1} z\\ 
-2 z^T I_{N_t} \otimes R_{x}^{-1} \diag \left(R_{x}^{-1}\right)^{-2} R_{x}^{-1} \frac{\partial R_{x}}{\partial \ell_{x_k}} R_{x}^{-1} z
\end{aligned},
\end{equation}
with
\begin{equation}\label{terLOOder}
\left(\frac{\partial R_{x,l}}{\partial \ell_{x_k}}\right)_{i,j} = \frac{\ell_{x_k}(x_{k,j}-x_{k,i})^2}{| x_{k,j}-x_{k,i} |}  h'_{\frac{5}{2}}\left(\frac{| x_{k,j}-x_{k,i} |}{\ell_{x_k}}\right)
\end{equation}
and
\begin{equation}\label{dermater52}
h'_{\frac{5}{2}}(x) =
-\frac{5}{3}x \left(1 + \sqrt{5}x\right)\exp\left(-\sqrt{5}x\right) .
\end{equation}
The method used to calculate the value of $\bm{\ell_x}$ is the Nelder-Mead method with only one starting point, because starting from more points will be more costly and the function $f_{\textbf{CV}}$ is close to quadratic consequently does not need multiple starting points.

\paragraph{Optimization without hypothesis (\ref{conditionmatcor}):} When hypothesis (\ref{conditionmatcor}) does not hold, a way must be found to calculate $\bm{\ell_x}$ without this assumption.
By a regularization of the matrix $R_t$ it is possible to calculate  $f_{\textbf{CV}}(\bm{\ell_x})$ and its derivative by \Cref{LOOcostfunctred,rtregularisation,eq::LOOerrnorm}.
However, the solution will be a regularized solution and not an exact solution.

There are different types of regularization that allow matrices to be inverted. Two methods have been investigated here.
The first one is standard (Tilkonov regularization): 
\begin{equation}\label{rtregularisation}
\widehat{R_t^{-1}} = \left( R_t^TR_t +\varepsilon^2 I_{N_t} \right)^{-1}R_t^T.
\end{equation}
The second one is:
\begin{equation}\label{regularisationaux}
\widehat{R_t^{-1}} = \left( R_t +\varepsilon I_{N_t} \right)^{-1}.
\end{equation}
It has the disadvantage of being more sensitive to $\varepsilon$ than the first one, which is why it will not be used. 

However, in the calculation of the determinant, this adjustment may have advantages.
\textcolor{black}{Denoting by} $\widehat{R_t^{-1}} = V \Sigma^{-1} U^T$ \textcolor{black}{the SVD of $R_t$,} with $\Sigma^{-1} = \diag{\frac{1}{\sigma_i+\varepsilon}}$ whereas the same decomposition gives for \Cref{rtregularisation} $\Sigma^{-1} = \diag{\frac{\sigma_i}{\sigma_i^2+\varepsilon^2}}$.
This is the reason why the two adjustments presented are not used in the same case.
Indeed $\frac{1}{\sigma_i+\varepsilon}$ is less efficient for the calculation of a determinant but more efficient for the calculation of the inverse of $\widehat{R_t}$.

\begin{equation}\label{LOOdercompl}
\begin{aligned}
\frac{\partial}{\partial \ell_{x_k}}f_{\textbf{CV}}(\bm{\ell_x}) =2 z^T \left(\widehat{R_t^{-1}} \diag \left(\widehat{R_t^{-1}}\right)^{-2}\widehat{R_t^{-1}} \right) \\
 \otimes R_{x}^{-1} \diag \left(R_{x}^{-1}\right)^{-2} \left( R_{x}^{-1} \frac{\partial R_{x}}{\partial \ell_{x_k}} R_{x}^{-1} \right) \diag \left(R_{x}^{-1}\right)^{-1} R_{x}^{-1} z\\ 
-2 z^T \left(\widehat{R_t^{-1}} \diag \left(\widehat{R_t^{-1}}\right)^{-2}\widehat{R_t^{-1}} \right) \otimes R_{x}^{-1} \diag \left(R_{x}^{-1}\right)^{-2} R_{x}^{-1} \frac{\partial R_{x}}{\partial \ell_{x_k}} R_{x}^{-1} z
\end{aligned}.
\end{equation}
\Cref{terLOOder} and \Cref{dermater52} are still valid.

This complete approach was compared to the LOO calculation using a loop.
However, the calculation time of Kronecker products is too long compared to the calculation of the simple error with a loop.
Moreover, the difference\textcolor{black}{s} in the errors of the different methods \textcolor{black}{are} negligible.
Thus this solution is only recommended when the calculation of \Cref{LOOcostfunctred,LOOdercompl} is optimized.

\begin{table}
\centering
\begin{tabular}{|l|c|c|c|}
  \hline
   & Loop LOO & Full LOO & Simplified LOO  \\
  \hline
  $\varepsilon_{\mathcal{Q^2}}$ & $7.06~10^{-4}$ & $8.72~10^{-4}$ & $8.14~10^{-4}$\\
  time & 18.41 s & 3.47 min & 0.17 s \\
  \hline
  
\end{tabular}
\caption{\label{bencLOOtech}Benchmark of the different LOO optimization techniques for estimating the $f(\bm{x},t) = \cos\left(4\pi (x_2+1)t\right)\sin\left(3\pi x_1t\right)$ function using the separable covariance method.
Loop LOO processes the error by computing the approximation, Full LOO processes the regularize\textcolor{black}{d} analytic expression and Simplified LOO processes the simplified one given by  \Cref{LOOcostsimp}.}
\end{table}

\Cref{bencLOOtech} shows that the gain in calculation time by the simplified method is significant even though the error difference is very small.
The extremely long time for the complete LOO is mainly due to an implementation of the Kronecker product that is not very effective in our implementation.

\section{Tensorized covariance of the orthogonal part}
\label{asec::TensOrt}
For $R_x$ we assume that $C_x$ is chosen in the Matérn-5/2 class of functions.
The function only depends on the correlation length vector $\bm{\ell_x}$.
The matrix $R_t$ is estimated as described in \Cref{subsec::tensRegcov} by the matrix $\widehat{R}_t$ given by:
$$\widehat{R_t} = \frac{1}{N_x}\left(\bm{Z}_{\text{obs}}^\perp - \bm{\hat{Z}}^\perp\right) R_x^{-1} \left(\bm{Z}_{\text{obs}}^\perp - \bm{\hat{Z}}^\perp\right)^T,$$
where $\bm{\hat{Z}}^\perp$ is the $N_t\times N_x$ matrix of empirical means $\hat{Z}_{u,i}^\perp=\frac{1}{N_x}\sum_{j=1}^{N_x}\left(\bm{Z}_\text{obs}^\perp\right)_{u,j}$, $\forall i=1, \ldots, N_x$ and $u=1,\ldots,N_t$. 
Its range is indeed in $S_N^\perp$.

The prediction mean is the sum of two terms, $\bm{Z}_\text{obs}^\perp R_x^{-1}r_x(\bm{x})$ which is $S_N^\perp$-valued and $B_\star u(\bm{x})$ also $S_N^\perp$-valued because:
\begin{equation}
    B_\star =\bm{Z_\text{obs}}^\perp R_x^{-1} F\left(F^T R_x^{-1}F\right)^{-1},
\end{equation}
with $F$ the $N_\text{F} \times M$ matrix $\left[f^T(\bm{x}^{(i)})\right]_{i=1,\hdots,N_\text{F}}$.
Consequently, we have:
\begin{equation}
    \label{eq::TenCovGammaN}
    Z_\text{H}^\perp(\bm{x},t_u) | \bm{\ell_x}, \bm{\Gamma}, N,\bm{Z}_\text{obs}^\perp \sim \mathcal{GP}(\mu_\star(\bm{x}) , R_\star(\bm{x},\bm{x}'))
\end{equation}
where the mean is given by (\ref{eq::meanTensCov}) and the covariance by (\ref{eq::varTensCov}) with $\bm{Z}_\text{obs}^\perp$ as the observed inputs.
LOO estimation of the vector of correlation lengths $\bm{\ell_x}$ given $\bm{\Gamma}$ and $N$ is carried out by the method presented in \Cref{asec::LOO}.

\section{\textcolor{black}{Expressions of some expectations and variances}}
\label{asec::Computation}
\subsection{Computation for uncorrolated Gaussian Processes}
\label{assec::MultiUnco}
Given $\bm{\Gamma}$,\\ $\left(A_{i,\text{H}}(\bm{x},t_u), A_{i,\text{L}}(\bm{x},t_u)\right)_{\substack{\bm{x}\in Q\\ u=1,\ldots,N_t}}$ are independent with respect to $i$.
This independence makes it possible to generate $N_t$ independent surrogate models, with mean and variance given by \cref{eq::FullUKMultiFi,eq::FullUKMultiVar}:

\begin{align}
    \label{eq::meanautomodeleq}
    \E\left[A_{i,\text{H}}(\bm{x})|\bm{\Gamma},\mathcal{A}=\alpha\right] = &
    h_{i,\text{H}}^T(\bm{x})\Sigma_{i,\text{H}}\nu_{i,\text{H}}
    + r_{i,\text{H}}^T(\bm{x})C_{i,\text{H}}^{-1}\left(\alpha_{i}^{\text{H}}-H_{i,\text{H}}\Sigma_{i,\text{H}}\nu_{i,\text{H}}\right), \\
\nonumber
        \V\left[A_{i,\text{H}}(\bm{x})|\bm{\Gamma},\mathcal{A}=\alpha\right]  =&
        \left(\hat{\rho}^2_{i,\text{L}}(\bm{x}) + \varepsilon_{i,\rho} (\bm{x})\right)\sigma_{\tilde{A}_{\text{L},i}}^2(\bm{x})\\
\nonumber        &
        + \frac{Q_{i,\text{H}}}{2(d_{i,\text{H}}-1)}(1-r_{i,\text{H}}^T(\bm{x})C_{i,\text{H}}^{-1}r_{i,\text{H}}(\bm{x}))\\
        &
 \hspace*{-1.0in}
        + \left(h_{i,\text{H}}^T(\bm{x})-r_{i,\text{H}}^T(\bm{x})C_{i,\text{H}}^{-1}H_{i,\text{H}}\right)\Sigma_{i,\text{H}}\left(h_{i,\text{H}}^T(\bm{x}) -r_{i,\text{H}}^T(\bm{x})C_{i,\text{H}}^{-1}H_{i,\text{H}}\right)^T .
    \label{eq::varautomodeleq}
\end{align}

\subsection{Variance for projection}
\label{assec:varproj}

The law of total variance gives :
\begin{equation}
        \V_\alpha\left[ Z_{\text{H}}(\bm{x},t_u) \right] =
        \V_\alpha\left[\E_\alpha\left[ Z_\text{H}(\bm{x},t_u)|\bm{\Gamma}\right]\right]
        + \E_\alpha\left[\V_\alpha\left[ Z_\text{H}(\bm{x},t_u)|\bm{\Gamma}\right]\right]
\end{equation}
The variance term can be expressed as follows :
\begin{equation}
    \label{eq::VarZparaVar}
    \begin{array}{c}
        \V_\alpha\left[\E_\alpha\left[ Z_\text{H}(\bm{x},t_u)|\bm{\Gamma}\right]\right] =
        \V_\alpha\left[ \sum_{i=1}^{N_t} \Gamma_i(t_u)\E_\alpha\left[A_{i,\text{H}}(\bm{x})|\bm{\Gamma}\right]\right] \\
        = \sum_{i=1}^{N_t} \V_\alpha\left[\Gamma_i (t_u)\E_\alpha \left[A_{i,\text{H}}(\bm{x})|\bm{\Gamma}\right]\right] \\
        + \sum_{i,j=1,i\neq j}^{N_t} \Cov_\alpha (\Gamma_i (t_u)\E_\alpha\left[A_{i,\text{H}}(\bm{x})|\bm{\Gamma} \right], \Gamma_j(t_u) \E_\alpha\left[A_{j,\text{H}}(\bm{x})|\bm{\Gamma} \right])
    \end{array},
\end{equation}
where $\E_\alpha \left[A_{i,\text{H}}(\bm{x})|\bm{\Gamma}\right]$ is given by \Cref{eq::meanautomodeleq}.
The expectation term can be expressed as :
\begin{equation}
    \label{eq::VarZparaMean}
    \hspace*{-0.2in}
    \begin{array}{c}
        \E_\alpha\left[\V_\alpha\left[ Z_\text{H}(\bm{x},t_u)|\bm{\Gamma}\right]\right] =\\
        \E_\alpha \left[ \sum_{i=1}^{N_t} \V_\alpha\left[A_{i,\text{H}}(\bm{x})\Gamma_i(t_u)| \bm{\Gamma}\right] + 
        \sum_{i,j=1,i\neq j}^{N_t}\Cov_\alpha\left(A_{i,\text{H}}(\bm{x})\Gamma_i(t_u),A_{j,\text{H}}(\bm{x})\Gamma_j(t_u)|\bm{\Gamma}\right)\right] \\
        = \sum_{i=1}^{N_t}\E_\alpha\left[\Gamma_i^2(t_u) \V_\alpha\left[A_{i,\text{H}}(\bm{x})|\bm{\Gamma}\right]\right] \\
        + \sum_{i,j=1,i\neq j}^{N_t} \E_\alpha\left[\Gamma_i(t_u)\Gamma_j(t_u)\Cov_\alpha\left(A_{i,\text{H}}(\bm{x}),A_{j,\text{H}}(\bm{x})|\bm{\Gamma}\right)\right]
    \end{array},
\end{equation}
where $\V_\alpha\left[A_{i,\text{H}}(\bm{x})\Gamma_i(t_u)| \bm{\Gamma}\right]$ is given in \Cref{eq::varautomodeleq} and \\ $ \Cov_\alpha(A_{i,\text{H}}(\bm{x}),A_{j,\text{H}}(\bm{x})|\bm{\Gamma})= 0$ 
if $i \neq j$.
Consequently:
\begin{equation}
    \begin{array}{c}
        \V_\alpha\left[ Z_{\text{H}}(\bm{x},t_u) \right] = 
        \sum_{i=1}^{N_t} \V_\alpha\left[\Gamma_i (t_u)\E_\alpha \left[A_{i,\text{H}}(\bm{x})|\bm{\Gamma}\right]\right] \\
        + \sum_{i,j=1,i\neq j}^{N_t} \Cov_\alpha (\Gamma_i (t_u)\E_\alpha\left[A_{i,\text{H}}(\bm{x})|\bm{\Gamma}\right], \Gamma_j(t_u) \E_\alpha\left[A_{j,\text{H}}(\bm{x})|\bm{\Gamma} \right])\\
        + \sum_{i=1}^{N_t}\E_\alpha\left[\Gamma_i^2(t_u) \V_\alpha\left[A_{i,\text{H}}(\bm{x})|\bm{\Gamma}\right]\right] \\
    \end{array}.
\end{equation}

\subsection{Variance for tensorisation of covariance and projection}
\label{assec:vartenso}

The theorem of the total variance gives us:
\begin{equation}
    \begin{array}{c}
        \V_{\bm{Z}_\text{obs}}\left[Z_{\text{H}}\left(\bm{x},t_u\right)|N,\bm{\ell_x}\right] =
        \V_{\bm{Z}_\text{obs}}\left[\E_{\bm{Z}_\text{obs}}\left[Z_{\text{H}}\left(\bm{x},t_u\right)|\bm{\Gamma},N,\bm{\ell_x}\right]|N,\bm{\ell_x}\right]\\
        + \E_{\bm{Z}_\text{obs}}\left[\V_{\bm{Z}_\text{obs}}\left[Z_{\text{H}}\left(\bm{x},t_u\right)|\bm{\Gamma},N,\bm{\ell_x}\right]|N,\bm{\ell_x}\right]
    \end{array}.
\end{equation}
The two terms of \cref{eq::NCalInterVar} are:
\begin{equation}
    \label{eq::VZobsEZobsfull}
    \begin{array}{c}
        \V_{\bm{Z}_\text{obs}}\left[\E_{\bm{Z}_\text{obs}}\left[Z_{\text{H}}\left(\bm{x},t_u\right)|\bm{\Gamma},N,\bm{\ell_x}\right]|N,\bm{\ell_x}\right]=\\ 
        \V_{\bm{Z}_\text{obs}}\left[\E_{\bm{Z}_\text{obs}}\left[Z_{\text{H}}^\perp\left(\bm{x},t_u\right)|\bm{\Gamma},N,\bm{\ell_x}\right] 
        + \E_{\bm{Z}_\text{obs}}\left[Z_{\text{H}}^\parallel\left(\bm{x},t_u\right)|\bm{\Gamma},N \right]|N \right]
    \end{array},
\end{equation}
and
\begin{equation}
    \label{eq::EZobsVZobsfull}
    \begin{array}{c}
     \E_{\bm{Z}_\text{obs}}\left[\V_{\bm{Z}_\text{obs}}\left[Z_{\text{H}}\left(\bm{x},t_u\right)|\bm{\Gamma},N,\bm{\ell_x}\right]|N,\bm{\ell_x}\right] =\\
     \E_{\bm{Z}_\text{obs}}\left[\V_{\bm{Z}_\text{obs}}\left[Z_{\text{H}}^\perp\left(\bm{x},t_u\right)|\bm{\Gamma}, N, \bm{\ell_x}\right]|N,\bm{\ell_x}\right]\\
      + \E_{\bm{Z}_\text{obs}}\left[2\text{Cov}_{\bm{Z}_\text{obs}}\left[ Z_{\text{H}}^\parallel\left(\bm{x},t_u\right), Z_{\text{H}}^\perp\left(\bm{x},t_u\right)|\bm{\Gamma}, N, \bm{\ell_x}\right]|N,\bm{\ell_x}\right] \\
      + \E_{\bm{Z}_\text{obs}}\left[\V_{\alpha}\left[Z_{\text{H}}^\parallel\left(\bm{x},t_u\right)|\bm{\Gamma},N\right]| N\right].
    \end{array}.
    \end{equation}

The uncorrelation of the $A_{i,\text{H}}(\bm{x},t_u)$ coefficients given $\bm{\Gamma}$ gives $\operatorname{Cov}_\alpha\left[A_{i,\text{H}}(\bm{x}),A_{j,\text{H}}(\bm{x})|\bm{\Gamma}\right] = 0$ for $i\neq j$ and $\operatorname{Cov}_{\bm{Z}_\text{obs}}\left[A_{i,\text{H}}(\bm{x})\Gamma_i(t_u),Z^\bot_\text{H}(\bm{x},t_u)|\bm{\Gamma}\right]=0$.
This leads us to simplify \Cref{eq::VZobsEZobsfull,eq::EZobsVZobsfull} into:
\begin{equation}
    \begin{array}{c}
        \V_{\bm{Z}_\text{obs}}\left[\E_{\bm{Z}_\text{obs}}\left[Z_\text{H}(\bm{x},t_u)|\bm{\Gamma}, N, \bm{\ell_x}\right]\right] =
        \V_{\bm{Z}_\text{obs}}\left[\E_{\alpha}\left[ Z_\text{H}^\perp(\bm{x},t_u)|\bm{\Gamma}, N, \bm{\ell_x}\right]| N, \bm{\ell_x}\right]\\
        + \sum_{i=1}^N \V_{\bm{Z}_\text{obs}}\left[\Gamma_i(t_u)\E_{\alpha}\left[A_{i,\text{H}}(\bm{x})|\bm{\Gamma}\right]\right] \\
        + \sum_{i,j=1; i\neq j}^N \text{Cov}_{\bm{Z}_\text{obs}}\left[\Gamma_i(t_u)\E_{\alpha}\left[A_{i,\text{H}}(\bm{x})|\bm{\Gamma}\right],\Gamma_j(t_u)\E_{\alpha}\left[A_{j,\text{H}}(\bm{x})|\bm{\Gamma}\right]\right]\\
        + 2\sum_{i=1}^N \text{Cov}_{\bm{Z}_\text{obs}}\left[\Gamma_i(t_u)\E_{\alpha}\left[A_{i,\text{H}}(\bm{x})|\bm{\Gamma}\right], \E_{\bm{Z}_\text{obs}}\left[Z_\text{H}^\perp(\bm{x},t_u)|\bm{\Gamma}, N, \bm{\ell_x}\right]| N, \bm{l_x}\right]
    \end{array},
\end{equation}
and 
\begin{equation}
    \begin{array}{c}
        \E_{\bm{Z}_\text{obs}}\left[\V_{\bm{Z}_\text{obs}}\left[Z_\text{H}\left(\bm{x},t_u\right)|\bm{\Gamma}, N, \bm{\ell_x}\right]\right] =
        \E_{\bm{Z}_\text{obs}}\left[\V_{\bm{Z}_\text{obs}}\left[Z_\text{H}^\perp\left(\bm{x},t_u\right)|\bm{\Gamma}, N, \bm{\ell_x}\right]| N, \bm{\ell_x}\right]\\
        + \sum_{i=1}^{N}\E_{\bm{Z}_\text{obs}}\left[\Gamma_i(t_u)^2\V_{\alpha}\left[A_{i,\text{H}}\left(\bm{x}\right)|\bm{\Gamma}\right]| N, \bm{\ell_x}\right]
    \end{array}.
\end{equation}

The full formula of the variance can be expressed as :
\begin{equation}
    \label{eq::VarFullortho}
    \begin{array}{c}
        \V_{\bm{Z}_\text{obs}}\left[Z_{\text{H}}\left(\bm{x},t_u\right)|N,\bm{\ell_x}\right] =
        \V_{\bm{Z}_\text{obs}}\left[\E_{\bm{Z}_\text{obs}}\left[ Z_\text{H}^\perp(\bm{x},t_u)|\bm{\Gamma},N,\bm{\ell_x}\right]|N,\bm{\ell_x}\right]\\
        + \E_{\bm{Z}_\text{obs}}\left[\V_{\alpha}\left[Z_\text{H}^\perp\left(\bm{x},t_u\right)|\bm{\Gamma}, N, \bm{\ell_x}\right]| N, \bm{\ell_x}\right]\\
        + \sum_{i=1}^N \V_{\bm{Z}_\text{obs}}\left[\Gamma_i(t_u)\E_{\alpha}\left[A_{i,\text{H}}(\bm{x})|\bm{\Gamma}\right]\right]\\ 
        + \sum_{i=1}^{N}\E_{\bm{Z}_\text{obs}}\left[\Gamma_i(t_u)^2\V_\alpha\left[A_i\left(\bm{x}\right)|\bm{\Gamma}\right]\right]\\
        + \sum_{i,j=1; i\neq j}^N \text{Cov}_{\bm{Z}_\text{obs}}\left[\Gamma_i(t_u)\E_\alpha\left[A_{i,\text{H}}(\bm{x})|\bm{\Gamma}\right],\Gamma_j(t_u)\E_\alpha\left[A_{j,\text{H}}(\bm{x})|\bm{\Gamma}\right]\right]\\
        + 2\sum_{i=1}^N \text{Cov}_{\bm{Z}_\text{obs}}\left[\Gamma_i(t_u)\E_\alpha\left[A_{i,\text{H}}(\bm{x})|\bm{\Gamma}\right], \E_{\bm{Z}_\text{obs}}\left[Z_\text{H}^\perp(\bm{x},t_u)|\bm{\Gamma}, N, \bm{\ell_x}\right]| N,  \bm{\ell_x}\right]\\
    \end{array}.
\end{equation}
where $ \V_\alpha{\left[A_{i,\text{H}}(\bm{x}) |\bm{\Gamma}\right]} $ is given by \Cref{eq::varautomodeleq} and $ \V_\alpha\left[Z_\text{H}^\perp(\bm{x})|\bm{\Gamma}\right] $ is given by \Cref{eq::varTensCov}.

\section*{Acknowledgments}
The author thanks Pr. Josselin Garnier and Dr. Claire Cannamela for their guidance and advice.

\bibliographystyle{siamplain.bst}
\bibliography{references}
\end{document}